\documentclass[11pt,reqno]{amsart}

\usepackage{amsmath,amsthm,amscd}
\usepackage{amssymb,txfonts}
\usepackage{euscript}
\usepackage[all]{xy}

\numberwithin{equation}{subsection}

\newcommand{\sgsp}{\renewcommand{\baselinestretch}{1}\tiny\normalsize}
\newcommand{\sqsp}{\renewcommand{\baselinestretch}{1.2}\tiny\normalsize}

\newcommand{\relphantom[1]}{\mathrel{\phantom{#1}}}

\newtheorem{thm}[subsection]{Theorem}

\newtheorem{cor}[subsection]{Corollary}

\newtheorem*{Thm}{Theorem}

\newcommand{\cat}[1]{{\EuScript #1}}
\newcommand{\cA}{\cat{A}}

\newcommand{\cF}{\cat{F}}

\newcommand{\bZ}{\mathbf{Z}}

\newcommand{\Lambdabar}{\overline{\Lambda}}

\newcommand{\Thetabar}{\overline{\Theta}}
\newcommand{\thetabar}{\overline{\theta}}

\newcommand{\Pibar}{\overline{\Pi}}

\newcommand{\Dbar}{\overline{D}}
\newcommand{\Rbar}{\overline{R}}

\newcommand{\sumprime}{\sideset{}{'}\sum}  

\newcommand{\CH}{\mathrm{CH}_{\mathrm{Hopf}}}
\newcommand{\CGS}{C_{\mathrm{GS}}}
\newcommand{\Cma}{C_{\mathrm{MA}}}
\newcommand{\Cmc}{C_{\mathrm{MC}}}
\newcommand{\Ccc}{C_{\mathrm{CC}}}
\newcommand{\Cca}{C_{\mathrm{CA}}}
\newcommand{\Ccb}{C_{\mathrm{CB}}}
\newcommand{\Cmb}{C_{\mathrm{MB}}}
\newcommand{\Hoch}{\mathrm{Hoch}}
\newcommand{\Hma}{H_{\mathrm{MA}}}
\newcommand{\Hmb}{H_{\mathrm{MB}}}
\newcommand{\Hcb}{H_{\mathrm{CB}}}
\newcommand{\Hmc}{H_{\mathrm{MC}}}
\newcommand{\Hca}{H_{\mathrm{CA}}}
\newcommand{\Hcc}{H_{\mathrm{CC}}}

\newcommand{\dma}{d_{\mathrm{MA}}}
\newcommand{\dI}{d_I}
\newcommand{\dII}{d_{II}}
\newcommand{\dIII}{d_{III}}

\DeclareMathOperator{\Id}{Id}
\DeclareMathOperator{\Hom}{Hom}
\DeclareMathOperator{\End}{End}
\DeclareMathOperator{\Aut}{Aut}
\DeclareMathOperator{\Der}{Der}
\DeclareMathOperator{\Coder}{Coder}
\DeclareMathOperator{\Bider}{Bider}
\DeclareMathOperator{\Ob}{Ob}

\begin{document}
\title{Deformation bicomplex of module-algebras}
\author{Donald Yau}

\begin{abstract}
The deformation bicomplex of a module-algebra over a bialgebra is constructed.  It is then applied to study algebraic deformations in which both the module structure and the algebra structure are deformed.  The cases of module-coalgebras, comodule-(co)algebras, and (co)module-bialgebras are also considered.
\end{abstract}


\email{dyau@math.ohio-state.edu}
\address{Department of Mathematics, The Ohio State University Newark, 1179 University Drive, Newark, OH 43055, USA}

\maketitle
\sqsp


\section{Introduction}
\label{sec:intro}

Let $H$ be a bialgebra.  An $H$-\emph{module-algebra} is an associative algebra $A$ that is also an $H$-module such that the multiplication map on $A$ becomes an $H$-module morphism.  This algebraic structure arises often in algebraic topology, quantum groups \cite[Chapter V.6]{kassel}, Lie and Hopf algebras theory \cite{d,mont,sweedler}, and group representations \cite[Chapter 3]{abe}.  For example, in algebraic topology, the complex cobordism $\mathrm{MU}^*(X)$ of a topological space $X$ is an $S$-module-algebra, where $S$ is the Landweber-Novikov algebra \cite{landweber,novikov} of stable cobordism operations.  Likewise, the singular mod $p$ cohomology $\mathrm{H}^*(X; \bZ/p)$ of a topological space $X$ is an $\cA_p$-module-algebra, where $\cA_p$ is the Steenrod algebra associated to the prime $p$ \cite{es,milnor}.  More examples of this form can be found in \cite{boardman}.

The purpose of this paper is twofold:
\begin{enumerate}
\item The deformation bicomplex $\Cma^{**}(A)$ for an $H$-module-algebra $A$ is constructed.
\item The deformation bicomplex is used to study algebraic deformations of $A$, where deformations are taken with respect to \emph{both} the $H$-module structure and the algebra structure on $A$.
\end{enumerate}
The deformation bicomplexes (respectively, tricomplexes) of module-coalgebras and comodule-(co)algebras (respectively, (co)module-bialgebras) are also constructed.

In \cite{yau}, the author studied algebraic deformations of module-algebras, in which only the $H$-module structure is deformed.  This paper generalizes \cite{yau}, which in turn is a generalization of \cite{yau1}.  The current deformation theory of (co)module-algebras (respectively, (co)module-coalgebras) also generalizes the classical deformation theory of associative algebras \cite{ger2} (respectively, coalgebras).  Moreover, deformation of a (co)module-bialgebra is a generalization of the Gerstenhaber-Schack deformation theory of a bialgebra \cite{gs1,gs2}.

In the deformation bicomplex $\Cma^{**}(A)$ of an $H$-module-algebra $A$, the $(p,q)$-entry is the module $\Hom(H^{\otimes q}, \Hom(A^{\otimes p},A))$.  The $0$th row $\Cma^{*,0}(A)$ coincides with the Hochschild cochain complex $\Hoch^*(A,A)$ of $A$ with coefficients in itself, which is the deformation complex of $A$ as an associative algebra \cite{ger2}.  In particular, it has a graded Lie bracket \cite{ger}.  It is not known, however, whether there is a graded Lie bracket on the whole bicomplex $\Cma^{**}(A)$.  In fact, since an $H$-module-algebra is not an algebra over a PROP, the general constructions of Markl \cite{markl} for PROPic algebras do not seem to apply.

Each higher row $\Cma^{*,q}(A) = \Hom(H^{\otimes q}, \Hom(A^{\otimes *},A))$ $(q \geq 1)$, though not a Hochschild cochain complex itself, is isomorphic to one, namely, $\Hoch^*(A, \Hom(H^{\otimes q},A))$.  Since the module $\Hom(H^{\otimes q},A)$ is an associative algebra, it induces an associative product on $\Hoch^*(A, \Hom(H^{\otimes q},A))$, making it into a differential graded associative algebra.  In particular, each row in $\Cma^{**}(A)$ is a differential graded associative algebra.  Each higher column $\Cma^{p,*}(A) = \Hoch^*(H, \Hom(A^{\otimes p},A))$ $(p \geq 1)$ is also a Hochschild cochain complex.  However, only the first column $\Cma^{1,*}$ is known to admit a non-trivial associative product, which is induced by the algebra structure on $\Hom(A,A)$.

Here is a summary of the various deformation bicomplexes/tricomplexes constructed in section \ref{sec:bicomplex} and sections \ref{sec:mc} - \ref{sec:cb}.

\begin{Thm}
The deformation bicomplex/tricomplex of $A$ is:
   \begin{enumerate}
   \item $\Cma^{**}(A) = \Hoch^*(H, \Hom(A^{\otimes *},A))$ if $A$ is an $H$-module-algebra;
   \item $\Cmc^{**}(A) = \Hoch^*(H, \Hom(A, A^{\otimes *}))$ if $A$ is an $H$-module-coalgebra;
   \item $\Cca^{**}(A) = \Hoch^*(A, H^{\otimes *} \otimes A)$ if $A$ is an $H$-comodule-algebra;
   \item $\Ccc^{**}(A) = \Hom(A, H^{\otimes *} \otimes A^{\otimes *})$ if $A$ is an $H$-comodule-coalgebra;
   \item $\Cmb^{***}(A) = \Hoch^*(H, \Hom(A^{\otimes *}, A^{\otimes *}))$ if $A$ is an $H$-module-bialgebra;
   \item $\Ccb^{***}(A) = \Hoch^*(A, H^{\otimes *} \otimes A^{\otimes *})$ if $A$ is an $H$-comodule-bialgebra.
   \end{enumerate}
\end{Thm}

It should be noted that our deformation bicomplex $\Cma^{**}(A)$ is different from the Hopf-Hochschild cochain complex $\CH^*(A,A)$ constructed by Kaygun \cite{kaygun}.  The construction of the cochain complex $\CH^*(A,A)$ is similar to the usual Hochschild cochain complex of $A$, but it also takes into account the $H$-linearity.  In particular, if $H$ is the ground field, then $\CH^*(A,A)$ coincides with the Hochschild cochain complex of $A$.  Moreover, $\CH^*(A,A)$ has the structure of a brace algebra with multiplication \cite{yau2}, which leads to a Gerstenhaber algebra structure on the Hopf-Hochschild cohomology of $A$.

\subsection{Organization}

The rest of this paper is organized as follows.

In section \ref{sec:prelim}, we fix notations that will be used throughout this paper.  The deformation bicomplex $\Cma^{**}(A)$ of an $H$-module-algebra $A$ is constructed in section \ref{sec:bicomplex}.  In section \ref{sec:cup}, the various $\cup$-products in $\Cma^{**}(A)$ are constructed.

Algebraic deformations of module-algebras are discussed in section \ref{sec:def}.  In particular, infinitesimals are properly identified with $2$-cocycles in the deformation bicomplex.  The vanishing of $H^2$ of the total complex of the deformation bicomplex implies rigidity, in the sense that every deformation is equivalent to the trivial one.

The constructions and arguments in the module-algebra case can be adapted to the cases of module-coalgebras, comodule-(co)algebras, and (co)module-bialgebras.  In sections \ref{sec:mc} - \ref{sec:cb}, we construct the deformation bicomplexes/tricomplexes and their $\cup$-products for these algebraic structures and state the corresponding deformation results.


\section{Preliminaries}
\label{sec:prelim}

The purposes of this preliminary section are to fix notations and to recall some basic facts about module-algebras and Hochschild cohomology.

\subsection{Notations}
\label{subsec:notations}

Fix a ground field $K$ once and for all.  Modules, linearity, $\Hom$, and $\otimes$ are all meant over $K$, unless otherwise specified.  Let $(H, \mu_H, \Delta_H)$ be a bialgebra with associative multiplication $\mu_H$ and coassociative comultiplication $\Delta_H$.

In a coassociative coalgebra $(C, \Delta)$, we use Sweedler's notation \cite{sweedler} for comultiplication:
   \[
   \Delta^p(x) \,=\, \sideset{}{_{(x)}}\sum x_{(1)} \otimes \cdots \otimes x_{(p+1)}
   \]
for $p \geq 1$.  The subscript in $\sum_{(x)}$ will sometimes be omitted.

Given an associative algebra $(A, \mu_A)$, a \emph{derivation on $A$} is a linear self-map $\varphi \in \End(A) = \Hom(A,A)$ such that
   \[
   \varphi(ab) \,=\, a\varphi(b) + \varphi(a)b
   \]
for all $a, b \in A$.  The module of all derivations on $A$ is denoted by $\Der(A)$, which is considered as a submodule of $\End(A)$.

\subsection{Module-algebra}
\label{subsec:module-algebra}

Let $(A, \mu_A)$ be an associative algebra.  Then $A$ is said to be an \emph{$H$-module-algebra} if and only if there exists an $H$-module structure $\lambda \in \Hom(H, \End(A))$ on $A$ such that $\mu_A$ becomes an $H$-module morphism.  In other words, $\lambda$ satisfies the following two conditions:
   \begin{equation}
   \label{eq:ma axioms}
   \begin{split}
   \lambda(xy) & \,=\, \lambda(x) \circ \lambda(y), \\
   \lambda(x)(ab) & \,=\, \sideset{}{_{(x)}}\sum \lambda(x_{(1)})(a) \cdot \lambda(x_{(2)})(b),
   \end{split}
   \end{equation}
for $x, y \in H$ and $a, b \in A$.

\subsection{Hochschild cohomology}
\label{subsec:Hoch}

Let $M$ be an $A$-bimodule.  The module of \emph{Hochschild $n$-cochains of $A$ with coefficients in $M$} \cite{hochschild} is defined to be
   \[
   \Hoch^n(A,M) \,\buildrel \text{def} \over=\, \Hom(A^{\otimes n}, M).
   \]
The coboundary map
   \[
   \delta^n_h \colon \Hoch^n(A,M) \to \Hoch^{n+1}(A,M)
   \]
is given by the alternating sum
   \[
   \delta^n_h \,=\, \sum_{i=0}^{n+1} \, (-1)^i \delta^n_h \lbrack i \rbrack,
   \]
where
   \begin{equation}
   \label{eq:Hoch diff}
   (\delta^n_h \lbrack i \rbrack)(\varphi) \,=\,
   \begin{cases}
   \alpha_l \circ (\Id_A \otimes \varphi) & \text{ if } i = 0, \\
   \varphi \circ (\Id_{A^{\otimes (i-1)}} \otimes \mu_A \otimes \Id_{A^{\otimes (n-i)}}) & \text{ if } 1 \leq i \leq n, \\
   \alpha_r \circ (\varphi \otimes \Id_A) & \text{ if } i = n+1
   \end{cases}
   \end{equation}
for $\varphi \in \Hoch^n(A,M)$.  Here $\alpha_l$ and $\alpha_r$ are the left and right actions of $A$ on $M$, respectively.  The $n$th cohomology module of $\Hoch^*(A,M)$ is denoted by $HH^n(A,M)$.

\subsection{DGA and Hochschild cup-product}
\label{subsec:Hoch cup}

A \emph{differential graded associative algebra} (DGA in what follows) $C = (C^*, d^*, \cup)$ consists of a cochain complex $(C^*, d^*)$ and an associative graded product on $C^*$ such that the Leibniz identity,
   \[
   d(x \cup y) \,=\, (dx) \cup y + (-1)^{\vert x \vert} x \cup (dy),
   \]
is satisfied for $x, y \in C^*$.

Suppose that $M$ is an $A$-bimodule and is an associative algebra itself.  Suppose in addition that the following three conditions are satisfied for all $a \in A$ and $m_1, m_2 \in M$:
   \begin{equation}
   \label{eq:3conditions}
   \begin{split}
   a(m_1 m_2) & \,=\, (am_1) m_2, \\
   (m_1 a)m_2 & \,=\, m_1(am_2), \\
   (m_1 m_2)a & \,=\, m_1(m_2a).
   \end{split}
   \end{equation}
Then it can be checked easily that the Hochschild cochain complex $\Hoch^*(A,M)$ becomes a DGA whose product, denoted by $\cup$, is given by
   \begin{equation}
   \label{eq:cup Hoch}
   (\varphi \cup \psi)(a_{1,r+s}) \,=\, \varphi(a_{1,r}) \cdot \psi(a_{r+1,r+s})
   \end{equation}
for $\varphi \in \Hoch^r(A,M)$, $\psi \in \Hoch^s(A,M)$, and $a_1, \ldots , a_{r+s} \in A$.  Here $a_{p,q}$ denotes the element $a_p \otimes \cdots \otimes a_q$ whenever $p \leq q$.  We will continue to use this shorthand throughout the rest of this paper.

The three conditions in \eqref{eq:3conditions} are satisfied, for example, when $f \colon A \to M$ is a morphism of algebras and $A$ acts on $M$ via $f$.


\section{Deformation bicomplex of module-algebras}
\label{sec:bicomplex}

From here on until the end of section \ref{sec:def}, $A$ will denote an $H$-module-algebra with $H$-module structure map $\lambda \in \Hom(H, \End(A))$.  The purpose of this section is to construct the deformation bicomplex $\Cma^{**}(A)$ of $A$.  Deformations, which will be discussed in section \ref{sec:def}, are taken with respect to both the $H$-module structure $\lambda$ and the multiplication structure $\mu_A$ on $A$.  As will be seen below, this bicomplex can be more explicitly denoted by $\Hoch^*(H,\Hom(A^{\otimes *},A))$.  Further structures of the deformation bicomplex will be discussed in section \ref{sec:cup}.

\subsection{$H$-bimodule structure on $\Hom(A^{\otimes n}, A)$}
\label{subsec:H-bimod}

When $A$ is an $H$-module-algebra, there is an $H$-bimodule structure on the module $\Hom(A^{\otimes n}, A)$ for each $n \geq 1$.  The left and right actions are given by
   \[
   \begin{split}
   (x \varphi)(a_{1,n}) &\,=\, \lambda(x)(\varphi(a_{1,n})), \\
   (\varphi x)(a_{1,n}) &\,=\, \sum_{(x)}\varphi(\lambda(x_{(1)})(a_1) \otimes \cdots \otimes \lambda(x_{(n)})(a_n))
   \end{split}
   \]
for $x \in H$, $\varphi \in \Hom(A^{\otimes n}, A)$, and $a_1, \ldots , a_n \in A$.

In particular, there is a Hochschild cochain complex $\Hoch^*(H, \Hom(A^{\otimes n}, A))$ of $H$ with coefficients in the $H$-bimodule $\Hom(A^{\otimes n}, A)$.

\subsection{Deformation bicomplex}
\label{subsec:bicomplex}

The deformation bicomplex of the $H$-module-algebra $A$ is the first quadrant, cohomological bicomplex
   \[
   \Cma^{**}(A) \,=\, \left\lbrace\Cma^{p,\,q}(A), d^{p,\,q}, (-1)^{p+1}b^{p,\,q}\right\rbrace
   \]
whose $(p,q)$-entry is
   \[
   \Cma^{p,\,q}(A) \,=\,
   \begin{cases}
   0 & \text{ if } p = 0, \\
   \Der(A) & \text{ if } (p,q) = (1,0), \\
   \Hom(H^{\otimes q}, \Hom(A^{\otimes p}, A)) & \text{ otherwise}.
   \end{cases}
   \]
The vertical and horizontal differentials in $\Cma^{**}(A)$ are denoted by
   \[
   \begin{split}
   (-1)^{p+1}b^{p,\,q} \colon & \Cma^{p,\,q}(A) \to \Cma^{p,\,q+1}(A), \\
   d^{p,\,q} \colon & \Cma^{p,\,q}(A) \to \Cma^{p+1,\, q}(A),
   \end{split}
   \]
respectively.  These differentials are defined as follows.  Note that
   \[
   \Cma^{p,\,q}(A) \,=\, \Hoch^q(H, \Hom(A^{\otimes p}, A)).
   \]
For $p \geq 1$, set
   \[
   b^{p,\,q} \,\buildrel \text{def} \over =\, \delta_h^q \colon
   \Hoch^q(H, \Hom(A^{\otimes p}, A)) \to \Hoch^{q+1}(H, \Hom(A^{\otimes p}, A)),
   \]
i.e., the Hochschild coboundary.  In particular, we have that
   \[
   b^{p,\,q} \,=\, \sum_{j=0}^{q+1} \, (-1)^j b^{p,\,q}\lbrack j \rbrack,
   \]
where
   \[
   b^{p,\,q}\lbrack j \rbrack \,=\, \delta_h^q\lbrack j \rbrack
   \]
as in \eqref{eq:Hoch diff}.  The only exception is that $b^{1,\, 0}$ is the restriction of $\delta_h^0$ to the submodule $\Der(A) \subseteq \Hom(A,A)$.  In what follows, such exceptions will be made for the $(1,0)$-entry automatically without further comments.

To define the horizontal differentials, note that
   \[
   \Cma^{p,\, 0}(A) \,=\, \Hom(A^{\otimes p}, A) \,=\, \Hoch^p(A,A).
   \]
In the $0$th row $\Cma^{*,\,0}(A)$, the differentials are defined by
   \[
   d^{p,\, 0} \,\buildrel \text{def} \over =\, \delta_h^p \colon
   \Hoch^p(A,A) \to \Hoch^{p+1}(A,A).
   \]
In particular, we have that
   \[
   d^{p,\,0} \,=\, \sum_{i=0}^{p+1} \, (-1)^i d^{p,\, 0}\lbrack i \rbrack,
   \]
where
   \[
   d^{p,\, 0}\lbrack i \rbrack \,=\, \delta_h^p \lbrack i \rbrack
   \]
as in \eqref{eq:Hoch diff}.

For $q \geq 1$, define the horizontal differential by
   \[
   d^{p,\,q}
   \,\buildrel \text{def} \over =\,
   \sum_{i=0}^{p+1} \, (-1)^i d^{p,\,q} \lbrack i \rbrack,
   \]
where
   \[
   (d^{p,\,q} \lbrack i \rbrack \varphi)(x_{1,q})(a_{1,p+1}) \,=\,
   \begin{cases}
   \sum \lambda(x_{1(1)} \cdots x_{q(1)})(a_1) \cdot \varphi(x_{1(2)} \otimes \cdots \otimes x_{q(2)})(a_{2,p+1}) & \text{ if } i = 0, \\
   \varphi(x_{1,q})(a_{1,i-1} \otimes (a_i a_{i+1}) \otimes a_{i+2,p+1}) & \text{ if } 1 \leq i \leq p, \\
   \sum \varphi(x_{1(1)} \otimes \cdots \otimes x_{q(1)})(a_{1,p}) \cdot \lambda(x_{1(2)} \cdots x_{q(2)})(a_{p+1}) & \text{ if } i = p+1.
   \end{cases}
   \]
Here $\varphi \in \Cma^{p,\,q}(A)$, $x_1, \ldots , x_q \in H$, and $a_1, \ldots , a_{p+1} \in A$.

So far we only know that the columns and the $0$th row in $\Cma^{**}(A)$ are cochain complexes.  To see that $\Cma^{**}(A)$ is a bicomplex, we need the following observations.

\begin{thm}
\label{thm:bicomplex ma}
\begin{enumerate}
\item For each $q \geq 1$ and $0 \leq k < l \leq p+2$, the equality
   \begin{equation}
   \label{eq:bicomplex ma d}
   d^{p+1,\, q} \lbrack l \rbrack \circ d^{p,\,q} \lbrack k \rbrack \,=\,
   d^{p+1,\, q} \lbrack k \rbrack \circ d^{p,\,q} \lbrack l-1 \rbrack
   \end{equation}
holds.
\item For $p \geq 1$, $q \geq 0$, $0 \leq i \leq p+1$, and $0 \leq j \leq q+1$, the equality
   \begin{equation}
   \label{eq:bicomplex ma db}
   d^{p,\, q+1} \lbrack i \rbrack \circ b^{p,\,q} \lbrack j \rbrack \,=\,
   b^{p+1,\, q} \lbrack j \rbrack \circ d^{p,\,q} \lbrack i \rbrack
   \end{equation}
holds.
\end{enumerate}
\end{thm}

\begin{proof}
Both statements can be proved by direct inspections on a case-by-case basis.  For instance, for \eqref{eq:bicomplex ma db} when $i = j = 0$, we have that
   \[
   \begin{split}
   (& b^{p+1,\, q} \lbrack 0 \rbrack \circ d^{p,\,q} \lbrack 0 \rbrack)(\varphi)(x_{1,q+1})(a_{1,p+1}) \\
   & =\, \lambda(x_1) \left( (d^{p,\,q} \lbrack 0 \rbrack)(\varphi)(x_{2,q+1})(a_{1,p+1}) \right) \\
   & =\, \lambda(x_1) \left( \sum \lambda(x_{2(1)} \cdots x_{q+1(1)})(a_1) \cdot \varphi(x_{2(2)} \otimes \cdots \otimes x_{q+1(2)})(a_{2,p+1})\right) \\
   & \buildrel (i) \over =\, \sum \lambda(x_{1(1)})\left(\lambda(x_{2(1)} \cdots x_{q+1(1)})(a_1)\right) \cdot \lambda(x_{1(2)})\left(\varphi(x_{2(2)} \otimes \cdots \otimes x_{q+1(2)})(a_{2,p+1})\right) \\
   & \buildrel (ii) \over =\, \sum \lambda(x_{1(1)} \cdots x_{q+1(1)})(a_1) \cdot \lambda(x_{1(2)})\left(\varphi(x_{2(2)} \otimes \cdots \otimes x_{q+1(2)})(a_{2,p+1})\right) \\
   & =\, \sum \lambda(x_{1(1)} \cdots x_{q+1(1)})(a_1) \cdot(b^{p,\,q} \lbrack 0 \rbrack \varphi)(x_{1(2)} \otimes \cdots \otimes x_{q+1(2)})(a_{2,p+1}) \\
   & =\, (d^{p,\, q+1} \lbrack 0 \rbrack \circ b^{p,\,q} \lbrack 0 \rbrack)(\varphi)(x_{1,q+1})(a_{1,p+1}).
   \end{split}
   \]
Here the equalities (i) and (ii) follow from the module-algebra axioms \eqref{eq:ma axioms}.

Similarly, if $i = 0$ and $1 \leq j \leq q$, then we have that
   \[
   \begin{split}
   (& b^{p+1,\, q} \lbrack j \rbrack \circ d^{p,\,q} \lbrack 0 \rbrack)(\varphi)(x_{1,q+1})(a_{1,p+1}) \\
   & =\, \sum \lambda(x_{1(1)} \cdots x_{q+1(1)})(a_1) \cdot \varphi\left(x_{1(2)} \otimes \cdots \otimes x_{j-1(2)} \otimes (x_{j(2)} x_{j+1(2)}) \otimes x_{j+2(2)} \otimes \cdots \otimes x_{q+1(2)}\right)(a_{2,p+1}) \\
   & =\, (d^{p,\, q+1} \lbrack 0 \rbrack \circ b^{p,\,q} \lbrack j \rbrack)(\varphi)(x_{1,q+1})(a_{1,p+1}).
   \end{split}
   \]
Note that only the bialgebra structure on $H$ is used in proving this condition.

Likewise, when $i = 0$ and $j = q+1$, we have that
   \[
   \begin{split}
   (& b^{p+1,\, q} \lbrack q+1 \rbrack \circ d^{p,\,q} \lbrack 0 \rbrack)(\varphi)(x_{1,q+1})(a_{1,p+1}) \\
   & =\, \sum \lambda(x_{1(1)} \cdots x_{q+1(1)})(a_1) \cdot \varphi(x_{1(2)} \otimes \cdots \otimes x_{q(2)})(\lambda(x_{q+1(2)})(a_2) \otimes \cdots \otimes \lambda(x_{q+1(p+1)})(a_{p+1})) \\
   & =\, (d^{p,\, q+1} \lbrack 0 \rbrack \circ b^{p,\,q} \lbrack q+1 \rbrack)(\varphi)(x_{1,q+1})(a_{1,p+1}).
   \end{split}
   \]
This proves \eqref{eq:bicomplex ma db} when $i = 0$.

If $1 \leq i \leq p$, then we have that
   \[
   \begin{split}
   (& d^{p,\, q+1} \lbrack i \rbrack \circ b^{p,\,q} \lbrack j \rbrack)(\varphi)(x_{1,q+1})(a_{1,p+1}) \\
   & =\,
   \begin{cases}
   \lambda(x_1) \left( \varphi(x_{2,q+1})(a_{1,i-1} \otimes (a_i a_{i+1}) \otimes a_{i+2,p+1})\right) & \text{ if } j = 0, \\
   \varphi(x_{1,j-1} \otimes (x_j x_{j+1}) \otimes x_{j+2,q+1})(a_{1,i-1} \otimes (a_i a_{i+1}) \otimes a_{i+2,p+1}) & \text{ if } 1 \leq j \leq q, \\
   \sum \varphi(x_{1,q})(\lambda(x_{q+1(1)})(a_1) \otimes \cdots \otimes \lambda(x_{q+1(i)})(a_i a_{i+1}) \otimes \cdots \otimes \lambda(x_{q+1(p)})(a_{p+1})) & \text{ if } j = q+1
   \end{cases}
   \\
   & =\, (b^{p+1,\, q} \lbrack j \rbrack \circ d^{p,\,q} \lbrack i \rbrack)(\varphi)(x_{1,q+1})(a_{1,p+1}).
   \end{split}
   \]
This proves \eqref{eq:bicomplex ma db} when $1 \leq i \leq p$.  This cases when $i = p + 1$ are similar to the cases when $i = 0$.

The condition \eqref{eq:bicomplex ma d} is proved by a similar analysis.
\end{proof}

From the condition \eqref{eq:bicomplex ma d}, it follows as usual that
   \[
   d^{p+1,\, q} \circ d^{p,\,q} \,=\, 0.
   \]
Therefore, each row $(\Cma^{*,\,q}(A), d^{*,\,q})$ is a cochain complex.  From the condition \eqref{eq:bicomplex ma db}, it follows that
   \[
   d^{p,\, q+1} \circ b^{p,\,q} \,=\,
   b^{p+1,\, q} \circ d^{p,\,q}.
   \]
In other words, each square in $\Cma^{**}(A)$ anti-commutes.  This leads to the following result.

\begin{cor}
\label{cor:Cma bicomplex}
The above definitions give a first quadrant, cohomological bicomplex $\Cma^{**}(A) = \left\lbrace \Cma^{p,\,q}, d^{p,\,q}, (-1)^{p+1}b^{p,\,q} \right\rbrace$.
\end{cor}

\begin{footnotesize}
   \[
   \UseTips
   \newdir{ >}{!/-5pt/\dir{>}}
   \xymatrix{
   \vdots & \vdots & \vdots & \vdots \\
   \Hoch^3(H, \End(A)) \ar[u]^{b^{1,3}} \ar[r]^{d^{1,3}} &
   \Hoch^3(H,\Hom(A^{\otimes 2},A)) \ar[u]^{-b^{2,3}} \ar[r]^{d^{2,3}} &
   \Hoch^3(H,\Hom(A^{\otimes 3},A)) \ar[u]^{b^{3,3}} \ar[r]^{d^{3,3}} &
   \Hoch^3(H,\Hom(A^{\otimes 4},A)) \ar[u]^{-b^{4,3}}  \\
   \Hoch^2(H, \End(A)) \ar[u]^{b^{1,2}} \ar[r]^{d^{1,2}} &
   \Hoch^2(H,\Hom(A^{\otimes 2},A)) \ar[u]^{-b^{2,2}} \ar[r]^{d^{2,2}} &
   \Hoch^2(H,\Hom(A^{\otimes 3},A)) \ar[u]^{b^{3,2}}  \ar[r]^{d^{3,2}} &
   \Hoch^2(H,\Hom(A^{\otimes 4},A)) \ar[u]^{-b^{4,2}} \\
   \Hoch^1(H, \End(A)) \ar[u]^{b^{1,1}} \ar[r]^{d^{1,1}} &
   \Hoch^1(H,\Hom(A^{\otimes 2,},A)) \ar[u]^{-b^{2,1}} \ar[r]^{d^{2,1}} &
   \Hoch^1(H,\Hom(A^{\otimes 3},A)) \ar[u]^{b^{3,1}} \ar[r]^{d^{3,1}} &
   \Hoch^1(H,\Hom(A^{\otimes 4},A)) \ar[u]^{-b^{4,1}} \\
   \Der(A) \ar[u]^{b^{1,0}} \ar[r]^{d^{1,0}=0} &
   \Hoch^2(A,A) \ar[u]^{-b^{2,0}} \ar[r]^{d^{2,0}} &
   \Hoch^3(A,A) \ar[u]^{b^{3,0}} \ar[r]^{d^{3,0}} &
   \Hoch^4(A,A) \ar[u]^{-b^{4,0}}
   }
   \]
\end{footnotesize}

The bicomplex $\Cma^{**}(A)$ is called the \emph{deformation bicomplex of $A$}.  Omitting the $0$th column, which is the $0$ cochain complex, the first four columns of $\Cma^{**}(A)$ appear as above.

Note that the deformation bicomplex $\Cma^{**}(A)$ contains the cochain complex $\cF^*(A)$ constructed in \cite{yau}.  In fact, $\cF^*(A)$ is, up to a shift in degree, the total complex of the sub-bicomplex of $\Cma^{**}(A)$ consisting of the first row $\Cma^{*,\,1}(A)$ and the first column $\Cma^{1,\,*}(A)$.

\subsection{Total complex and cohomology}
\label{subsec:total}

Denote by $(\Cma^*(A), \dma^*)$ the total complex of the deformation bicomplex $\Cma^{**}(A)$.  In particular:
   \[
   \Cma^n(A) \,=\,
   \begin{cases}
   0 & \text{ if } n = 0, \\
   \Der(A) & \text{ if } n = 1, \\
   \bigoplus_{i=1}^n \Hoch^{n-i}(H, \Hom(A^{\otimes i},A)) & \text{ if } n \geq 2.
   \end{cases}
   \]
The cochain complex $(\Cma^*(A), \dma^*)$ is called the \emph{deformation complex of $A$}.

Define the cohomology module
   \[
   \Hma^n(A) \, \buildrel \text{def} \over =\,
   H^n(\Cma^*(A), \dma^*).
   \]
The deformation complex and its cohomology modules will be used in section \ref{sec:def} to study algebraic deformations of $A$.


\section{Cup-product}
\label{sec:cup}

In this section, we observe that each row in the deformation bicomplex $\Cma^{**}(A)$ is a differential graded associative algebra.  Moreover, each row $\Cma^{*,\,q}(A)$ is canonically isomorphic to a Hochschild cochain complex.  Under this isomorphism, the product in each row corresponds to a Hochschild $\cup$-product.

\subsection{Cup-product in $\Cma^{*,\,q}(A)$}
\label{subsec:cup ma}

The Hochschild $\cup$-product in the $0$th row $\Cma^{*,\,0}(A) = \Hoch^*(A,A)$ generalizes to the higher rows $\Cma^{*,\,q}(A)$ $(q \geq 1)$.  In fact, one can define a product
   \[
   - \cup - \colon \Hom(H^{\otimes q}, \Hom(A^{\otimes r},A)) \otimes \Hom(H^{\otimes q}, \Hom(A^{\otimes s},A))
   \to \Hom(H^{\otimes q}, \Hom(A^{\otimes (r+s)},A))
   \]
by setting
   \begin{equation}
   \label{eq:cupma}
   (\varphi \cup \psi)(x_{1,q})(a_{1,r+s})
   \, \buildrel \text{def} \over = \,
   \sum \varphi(x_{1(1)} \otimes \cdots \otimes x_{q(1)})(a_{1,r}) \cdot \psi(x_{1(2)} \otimes \cdots \otimes x_{q(2)})(a_{r+1,r+s})
   \end{equation}
for $\varphi \in \Hom(H^{\otimes q}, \Hom(A^{\otimes r},A))$, $\psi \in \Hom(H^{\otimes q}, \Hom(A^{\otimes s},A))$, $x_i \in H$, and $a_j \in A$.

\begin{thm}
\label{thm:cup ma}
The product defined in \eqref{eq:cupma} is associative and satisfies the Leibniz identity,
   \begin{equation}
   \label{eq:leibniz ma}
   d^{r+s,\, q}(\varphi \cup \psi) \,=\, (d^{r,\, q}\varphi) \cup \psi + (-1)^r \varphi \cup (d^{s,\, q}\psi).
   \end{equation}
In particular, $\left(\Cma^{*,\,q}(A), d^{*,\,q}, \cup\right)$ is a DGA.
\end{thm}

\begin{proof}
The associativity of the $\cup$-product is a consequence of the coassociativity of the comultiplication $\Delta_H$ in $H$ and the associativity of the multiplication $\mu_A$ in $A$.  The Leibniz identity can be checked by a simple inspection of both sides of \eqref{eq:leibniz ma} when applied to $x_{1,q} \in H^{\otimes q}$ and then to $a_{1,r+s+1} \in A^{\otimes (r+s+1)}$.
\end{proof}

The Leibniz identity implies, as usual, that the $\cup$-product descends to cohomology, which leads to the following consequence.

\begin{cor}
\label{cor:cup ma}
For each $q \geq 0$, $\left(\bigoplus_{n \geq 1} H^n(\Cma^{*,\,q}(A), d^{*,\,q}), \cup\right)$ is a graded associative algebra.
\end{cor}

There is a more conceptual way to obtain the $\cup$-product above by realizing it as a Hochschild $\cup$-product, which we do next.

\subsection{Alternative description of $\Cma^{*,\,q}(A)$}
\label{subsec:alternate}

Fix an integer $q \geq 1$.  Using the second module-algebra axiom in \eqref{eq:ma axioms}, one obtains an $A$-bimodule structure on the module $\Hom(H^{\otimes q},A)$ via the actions
   \[
   \begin{split}
   (a f)(x_{1,q}) & \,=\, \sum \lambda(x_{1(1)} \cdots x_{q(1)})(a) \cdot f(x_{1(2)} \otimes \cdots \otimes x_{q(2)}), \\
   (f a)(x_{1,q}) & \,=\, \sum f(x_{1(1)} \otimes \cdots \otimes x_{q(1)}) \cdot \lambda(x_{1(2)} \cdots x_{q(2)})(a)
   \end{split}
   \]
for $a \in A$, $f \in \Hom(H^{\otimes q},A)$, and $x_1, \ldots , x_q \in H$.

Moreover, the module $\Hom(H^{\otimes q},A)$ is an associative algebra via the product
   \[
   (f \cdot g)(x_{1,q}) \,=\, \sum f(x_{1(1)} \otimes \cdots \otimes x_{q(1)}) \cdot g(x_{1(2)} \otimes \cdots \otimes x_{q(2)})
   \]
for $f, g \in \Hom(H^{\otimes q},A)$ and $x_i \in H$.  The three conditions in \eqref{eq:3conditions} can be checked easily in this case.  It follows as in section \ref{subsec:Hoch cup} that the Hochschild cochain complex $\Hoch^*(A, \Hom(H^{\otimes q},A))$ is a DGA whose associative $\cup$-product is given by
   \[
   (\varphi \cup \psi)(a_{1,r+s})(x_{1,q}) \,=\, \sum \varphi(a_{1,r})(x_{1(1)} \otimes \cdots \otimes x_{q(1)}) \cdot \psi(a_{r+1,r+s})(x_{1(2)} \otimes \cdots \otimes x_{q(2)})
   \]
for $\varphi \in \Hoch^r(A, \Hom(H^{\otimes q},A))$, $\psi \in \Hoch^s(A, \Hom(H^{\otimes q},A))$, $a_i \in A$, and $x_j \in H$.

\begin{thm}
\label{thm:C*q}
There is a canonical isomorphism
   \[
   \zeta \colon
   \left(\Hoch^p(A, \Hom(H^{\otimes q},A)), \delta_h^p, \cup\right)_{p \geq 1}
   \,\xrightarrow{\cong}\,
   \left(\Cma^{*,\,q}(A), d^{*,\,q}, \cup\right)
   \]
of DGAs defined by
   \begin{equation}
   \label{eq:zeta}
   (\zeta \varphi)(x_{1,q})(a_{1,p}) \buildrel \mathrm{def} \over = \varphi(a_{1,p})(x_{1,q})
   \end{equation}
for $\varphi \in \Hoch^p(A, \Hom(H^{\otimes q},A))$, $x_i \in H$, and $a_j \in A$.
\end{thm}

\begin{proof}
It is clear that the map $\zeta$ defined in \eqref{eq:zeta} is a linear isomorphism.  In fact, it is just the usual $\Hom-\otimes$ adjunction applied twice.  Direct inspections then show that, under the isomorphism $\zeta$, $\delta_h^*$ corresponds to $d^{*,\,q}$ and the $\cup$-product in $\Hoch^*(A, \Hom(H^{\otimes q},A))$ corresponds to the one in $\Cma^{*,\,q}(A)$.
\end{proof}

Passing to cohomology, this leads to the following result.

\begin{cor}
\label{cor:cup ma coh}
The isomorphism $\zeta$ \eqref{eq:zeta} induces an isomorphism
   \[
   \zeta \colon
   \left(\bigoplus_{n \geq 1} HH^n(A, \Hom(H^{\otimes q}, A)), \cup\right)
   \, \xrightarrow{\cong} \,
   \left(\bigoplus_{n \geq 1} H^n(\Cma^{*,\,q}(A), d^{*,\,q}), \cup\right)
   \]
of graded associative algebras.
\end{cor}

Note that the $\cup$-product in the first row $\Cma^{*,\,1}(A)$ coincides with the one in $\cF^*_1(A)$ constructed in \cite{yau}.


\section{Deformations of module-algebras}
\label{sec:def}

In this section, we show that the deformation complex $\Cma^*(A)$ is actually the cochain complex that controls the deformations, in the sense of Gerstenhaber \cite{ger2}, of an $H$-module-algebra $A$, in which both the $H$-module structure $\lambda \in \Hom(H, \End(A))$ and the multiplication $\mu_A \in \Hom(A^{\otimes 2},A)$ on $A$ are deformed.  This generalizes the treatment in \cite{yau}, in which only the $H$-module structure is deformed.  However, once the correct definitions are made, the arguments here are similar to those in \cite{yau}.

\subsection{Deformation}
\label{subsec:deformation}

A \emph{deformation of $A$} is a power series
   \[
   \Theta \,=\, \sum_{n\geq0} \theta_n t^n,
   \]
with $\theta_0 = (\lambda, \mu_A) \in \Cma^2(A)$ and each $\theta_n = (\lambda_n, \pi_n) \in \Cma^2(A)$, satisfying the following three conditions.  Writing $\Lambda = \sum_{n\geq 0} \lambda_n t^n$  $(\lambda_0 = \lambda)$ and $\Pi = \sum_{n\geq 0} \pi_n t^n$ $(\pi_0 = \mu_A)$, the three conditions are:
   \begin{subequations}
   \label{eq:def ma axioms}
   \begin{align}
   \Lambda(xy) & \,=\, \Lambda(x) \circ \Lambda(y), \label{eq:def ma 1} \\
   \Lambda(x)\left(\Pi(a,b)\right) & \,=\, \sideset{}{_{(x)}}\sum \Pi\left(\Lambda(x_{(1)})(a), \Lambda(x_{(2)})(b)\right), \label{eq:def ma 2} \\
   \Pi\left(\Pi(a,b),c\right) & \,=\, \Pi\left(a, \Pi(b,c)\right) \label{eq:def ma 3}
   \end{align}
   \end{subequations}
for $x, y \in H$ and $a, b, c \in A$.  Such a deformation will also be denoted by $\Theta = (\Lambda, \Pi)$.

The linear coefficient $\theta_1 = (\lambda_1, \pi_1)$ of a deformation $\Theta$ is called the \emph{infinitesimal}.  The \emph{trivial deformation} is the deformation $\Theta = \theta_0 = (\lambda, \mu_A)$.

\subsection{Equivalence}
\label{subseq:equivalence}

A \emph{formal automorphism of $A$} is a power series
   \[
   \Phi \,=\, \sum_{n\geq 0} \phi_n t^n,
   \]
in which $\phi_0 = \Id_A$ and each $\phi_n \in \End(A)$, such that the first non-zero $\phi_n$ $(n \geq 1)$ is a derivation on $A$.  Note that $\phi_1$ is necessarily a derivation on $A$.

Two deformations $\Theta = (\Lambda, \Pi)$ and $\Thetabar = (\Lambdabar, \Pibar)$ are said to be \emph{equivalent} if and only if there exists a formal automorphism $\Phi$ such that the following two conditions are satisfied:
   \begin{subequations}
   \label{eq:equiv ma}
   \begin{align}
   \Pibar & \,=\, \Phi^{-1} \circ \Pi \circ \Phi^{\otimes 2}, \label{eq:equiv ma 1} \\
   \Lambdabar & \,=\, \Phi^{-1} \Lambda \Phi. \label{eq:equiv ma 2}
   \end{align}
   \end{subequations}
On the right-hand side of \eqref{eq:equiv ma 2}, we use the interpretation
   \[
   (\phi_i \lambda_j \phi_k)(x)
   \,=\,
   \phi_i \circ \lambda_j(x) \circ \phi_k
   \]
for $x \in H$.  In the above situation, we write $\Thetabar = \Phi^{-1} \Theta \Phi$.  This defines an equivalence relation.

Given a deformation $\Theta = (\Lambda, \Pi)$ and a formal automorphism $\Phi$, one can define an equivalent deformation $\Thetabar = \Phi^{-1} \Theta \Phi$ using \eqref{eq:equiv ma 1} and \eqref{eq:equiv ma 2}.

The $H$-module-algebra $A$ is said to be \emph{rigid} if every deformation of $A$ is equivalent to the trivial deformation.

\begin{thm}
\label{thm:def ma}
Let $A$ be an $H$-module-algebra.  Then the following statements hold.
\begin{enumerate}
\item The infinitesimal $\theta_1$ of any deformation $\Theta$ of $A$ is a $2$-cocycle in the deformation complex $\Cma^2(A)$ whose cohomology class is determined by the equivalence class of $\Theta$.
\item If $\Hma^2(A) = 0$, then $A$ is rigid.
\end{enumerate}
\end{thm}

\begin{proof}
The deformation arguments in \cite{ger2} and \cite{yau} can be adapted to the present situation.  For example, the deformation axioms \eqref{eq:def ma axioms} can be rewritten as:
   \begin{subequations}
   \label{eq:def ma axioms'}
   \begin{align}
   \lambda_n(xy) & \,=\, \sum_{i+j=n} \lambda_i(x) \circ \lambda_j(y), \label{eq:def ma 1'} \\
   \sum_{i+j=n} \lambda_i(x)\left(\pi_j(a,b)\right) &\,=\, \sum_{(x)} \sum_{i+j+k=n} \pi_i\left(\lambda_j(x_{(1)})(a), \lambda_k(x_{(2)})(b)\right), \label{eq:def ma2'} \\
   \sum_{i+j=n}\pi_i\left(\pi_j(a,b),c\right) &\,=\, \sum_{i+j=n} \pi_i\left(a, \pi_j(b,c)\right) \label{eq:def ma 3'}
   \end{align}
   \end{subequations}
for $n \geq 1$, $x, y \in H$, and $a, b, c \in A$.  When $n = 1$, these three conditions state that
   \[
   \begin{split}
   b^{1,\, 1}\lambda_1 & \,=\, 0, \\
   d^{1,\, 1}\lambda_1 - b^{2,\, 0}\pi_1 & \,=\, 0, \\
   d^{2,\, 0}\pi_1 & \,=\, 0.
   \end{split}
   \]
These three statements together state that $\theta_1 = (\lambda_1, \pi_1) \in \Cma^2(A)$ is a $2$-cocycle.  Moreover, if $\Thetabar = \Phi^{-1} \Theta \Phi$ for some formal automorphism $\Phi$, then the condition on the linear coefficients can be restated as
   \[
   \thetabar_1 - \theta_1 \,=\, \dma^1\phi_1,
   \]
which is a $2$-coboundary in $\Cma^2(A)$.  This proves statement (1).  Statement (2) is proved similarly by adapting Proposition 3.5 in \cite{yau}.
\end{proof}


\section{Deformation bicomplex of module-coalgebras}
\label{sec:mc}

In this section, the deformation bicomplex $\Cmc^{**}(A) = \Hoch^*(H,\Hom(A,A^{\otimes *}))$ for an $H$-module-coalgebra $A$ is constructed, in which both the $H$-module structure and the coalgebra structure on $A$ are deformed.  The corresponding deformation results are then listed.  The proofs can be adapted from the module-algebra case.

\subsection{Module-coalgebra}
\label{subsec:mc}

Throughout this section, let $A = (A, \Delta_A)$ be a coassociative coalgebra.  A \emph{coderivation on $A$} is a linear self-map $\varphi \in \Hom(A,A)$ such that
   \[
   \Delta_A \circ \varphi \,=\, (\Id_A \otimes \varphi + \varphi \otimes \Id_A) \circ \Delta_A.
   \]
The set of coderivations on $A$ is denoted by $\Coder(A)$, which is considered as a submodule of $\Hom(A,A)$.

An \emph{$H$-module-coalgebra} structure on $A$ consists of an $H$-module structure $\lambda \in \Hom(H, \End(A))$ on $A$ such that the map $\Delta_A \colon A \to A \otimes A$ becomes an $H$-module morphism, i.e.,
   \[
   \Delta_A(\lambda(x)(a)) \,=\, \sum_{(a)(x)} \lambda(x_{(1)})(a_{(1)}) \otimes \lambda(x_{(2)})(a_{(2)})
   \]
for all $x \in H$ and $a \in A$.  For the rest of this section, $A$ will denote an $H$-module-coalgebra with $H$-module structure map $\lambda$.

\subsection{Hochschild coalgebra cohomology}
\label{subsec:coHoch}

The deformation bicomplex of $A$ uses the coalgebra version of Hochschild cohomology \cite{jonah,pw}, which we now recall.

Let $M$ be an $A$-bicomodule with left $A$-coaction $\psi_l$ and right $A$-coaction $\psi_r$.  Define the module of \emph{Hochschild coalgebra $n$-cochains of $A$ with coefficients in $M$} as
   \[
   \Hoch^n_c(M,A) \,=\,
   \begin{cases}
   0 & \text{ if } n = 0, \\
   \Hom(M, A^{\otimes n}) & \text{ if } n \geq 1.
   \end{cases}
   \]
The coboundary is defined by
   \[
   \delta_c \sigma
   \,=\,
   (\Id_A \otimes \sigma) \circ \psi_l  \,+\,
   \sum_{i=1}^n \, (-1)^i \left(\Id_{A^{\otimes(i-1)}} \otimes \Delta_A \otimes \Id_{A^{\otimes(n-i)}}\right) \circ \sigma +
   (-1)^{n+1} (\sigma \otimes \Id_A) \circ \psi_r
   \]
for $\sigma \in \Hoch^n_c(M,A)$.

Suppose, in addition, that $M = (M, \Delta_M)$ is a coassociative coalgebra such that the following three conditions, dual to \eqref{eq:3conditions}, are satisfied:
   \begin{equation}
   \label{eq:3cond}
   \begin{split}
   (\Id_A \otimes \Delta_M) \circ \psi_l & \,=\,  (\psi_l \otimes \Id_M) \circ \Delta_M, \\
   (\Delta_M \otimes \Id_A) \circ \psi_r & \,=\, (\Id_M \otimes \psi_r) \circ \Delta_M, \\
   (\psi_r \otimes \Id_M) \circ \Delta_M & \,=\, (\Id_M \otimes \alpha_l) \circ \Delta_M.
   \end{split}
   \end{equation}
Then the Hochschild coalgebra cochain complex $\Hoch_c^*(M,A)$ becomes a DGA with the product
   \begin{equation}
   \label{eq:Hoch c cup}
   f \cup g \,=\, (f \otimes g) \circ \Delta_M
   \end{equation}
for $f, g \in \Hoch_c^*(M,A)$.

\subsection{$H$-bimodule structure on $\Hom(A,A^{\otimes n})$}
\label{subsec:Hom(A,An)}

For $n \geq 1$, there is an $H$-bimodule structure on $\Hom(A,A^{\otimes n})$ defined as follows.  The left and right $H$-actions on $\Hom(A,A^{\otimes n})$ are given by
   \[
   \begin{split}
   (x\varphi)(a)
   &\,=\, \sum \lambda(x_{(1)})(\varphi(a)^1) \otimes \cdots \otimes \lambda(x_{(n)})(\varphi(a)^n), \\
   (\varphi x)(a)
   &\,=\, \varphi(\lambda(x)(a))
   \end{split}
   \]
for $x \in H$, $\varphi \in \Hom(A, A^{\otimes n})$, and $a \in A$.  In the left $H$-action, we use the notation
   \[
   \varphi(a) \,=\,
   \sum \varphi(a)^1 \otimes \cdots \otimes \varphi(a)^n \in A^{\otimes n}.
   \]
In particular, we can consider the Hochschild cochain complex $\Hoch^*(H,\Hom(A,A^{\otimes n}))$.

\subsection{Module-coalgebra deformation bicomplex}
\label{subsec:bicomplex mc}

The \emph{deformation bicomplex} of the $H$-module-coalgebra $A$ is the first quadrant, cohomological bicomplex
   \[
   \Cmc^{**}(A) \,=\, \left\lbrace \Cmc^{p,\, q}(A), d^{p,\, q}, (-1)^{p+1}b^{p,\, q}\right\rbrace
   \]
with
   \[
   \Cmc^{p,\, q}(A) \,=\,
   \begin{cases}
   0 & \text{ if } p = 0, \\
   \Coder(A) & \text{ if } (p,q) = (1,0), \\
   \Hom(H^{\otimes q}, \Hom(A, A^{\otimes p})) & \text{ otherwise}.
   \end{cases}
   \]
For each $p \geq 1$, set
   \[
   b^{p,\, q} \,\buildrel \text{def} \over=\, \delta_h^q \colon \Hoch^q(H,\Hom(A,A^{\otimes p})) \to \Hoch^{q+1}(H,\Hom(A,A^{\otimes p})).
   \]

In the $0$th row $\Cmc^{*,\,0}(A)$, define the horizontal differential as the Hochschild coalgebra coboundary,
   \[
   d^{*,\,0} \, \buildrel \text{def} \over=\, \delta_c^* \colon \Hoch_c^*(A,A) \to \Hoch_c^{*+1}(A,A).
   \]
For $q \geq 1$, the horizontal differential
   \[
   d^{p,\, q} \colon \Cmc^{p,\, q}(A) \to \Cmc^{p+1,\, q}(A)
   \]
is defined as the alternating sum
   \[
   d^{p,\, q} \,=\, \sum_{i=0}^{p+1} \, (-1)^i d^{p,\, q} \lbrack i \rbrack,
   \]
where
   \[
   (d^{p,\, q} \lbrack i \rbrack \varphi)(x_{1,q})(a)
   \,=\,
   \begin{cases}
   \sum \lambda(x_{1(1)} \cdots x_{q(1)})(a_{(1)}) \otimes \varphi(x_{1(2)} \otimes \cdots \otimes x_{q(2)})(a_{(2)}) & \text{ if } i = 0, \\
   (\Id_{A^{\otimes(i-1)}} \otimes \Delta_A \otimes \Id_{A^{\otimes (p-i)}})(\varphi(x_{1,q})(a)) & \text{ if } 1 \leq i \leq p, \\
   \sum \varphi(x_{1(1)} \otimes \cdots \otimes x_{q(1)})(a_{(1)}) \otimes \lambda(x_{1(2)} \cdots x_{q(2)})(a_{(2)}) & \text{ if } i = p+1.
   \end{cases}
   \]
Here $\varphi \in \Cmc^{p,\, q}(A)$, $x_j \in H$, and $a \in A$.

Theorem \ref{thm:bicomplex ma} still holds in the module-coalgebra case with essentially the same proof. It follows that $\Cmc^{**}(A)$ is indeed a bicomplex.  Denote its total complex by $\Cmc^*(A)$, which is called the \emph{deformation complex} of $A$.  The $n$th cohomology module of $\Cmc^*(A)$ is denoted by $\Hmc^n(A)$.

Note that the cochain complex $\cF^*_{mc}(A)$ constructed in \cite{yau} is the total complex of the sub-bicomplex of $\Cmc^{**}(A)$ consisting of the first column $\Cmc^{1,\,*}(A)$ and the first row $\Cmc^{*,\,1}(A)$.

\subsection{Cup-product}
\label{subsec:cup mc}

As in the module-algebra case, each row $(\Cmc^{*,\,q}(A), d^{*,\,q})$ is a DGA whose product is defined by
   \[
   (\varphi \cup \psi)(x_{1,q})(a) \,\buildrel \text{def} \over=\,
   \sum \varphi(x_{1(1)} \otimes \cdots \otimes x_{q(1)})(a_{(1)}) \otimes \psi(x_{1(2)} \otimes \cdots \otimes x_{q(2)})(a_{(2)})
   \]
for $\varphi \in \Cmc^{r,\,q}(A)$ and $\psi \in \Cmc^{s,\,q}(A)$.

\subsection{Module-coalgebra deformation}
\label{subsec:deformation mc}

A \emph{deformation} of the $H$-module-coalgebra $A$ is a power series $\Theta = \sum_{n \geq 0}\theta_n t^n$ with $\theta_0 = (\lambda, \Delta_A)$ and each $\theta_n = (\lambda_n, \Delta_n) \in \Cmc^2(A)$ such that the following three conditions are satisfied.  Writing $\Lambda = \sum_{n \geq 0} \lambda_n t^n$ $(\lambda_0 = \lambda)$ and $D = \sum_{n \geq 0}\Delta_n t^n$ $(\Delta_0 = \Delta_A)$, the three conditions are:
   \begin{equation}
   \label{eq:def mc}
   \begin{split}
   \Lambda(xy) & \,=\, \Lambda(x) \circ \Lambda(y), \\
   D \circ \Lambda(x) & \,=\, \left\lbrace \sideset{}{_{(x)}}\sum \Lambda(x_{(1)}) \otimes \Lambda(x_{(2)}) \right\rbrace \circ D,  \\
   (D \otimes \Id_A) \circ D & \,=\, (\Id_A \otimes D) \circ D
   \end{split}
   \end{equation}
for $x, y \in H$.  Such a deformation is also denoted by $\Theta = (\Lambda, D)$.  The linear coefficient $\theta_1 \in \Cmc^2(A)$ is called the \emph{infinitesimal}.  The \emph{trivial deformation} is the deformation $\Theta = \theta_0 = (\lambda, \Delta_A)$.

A \emph{formal automorphism} of $A$ is a power series $\Phi = \sum_{n \geq 0} \phi_n t^n$ with $\phi_0 = \Id_A$ and each $\phi_n \in \End(A)$ such that the first non-zero $\phi_n$ $(n \geq 1)$ is a coderivation on $A$.

Two deformations $\Theta = (\Lambda, D)$ and $\Thetabar = (\Lambdabar, \Dbar)$ are said to be \emph{equivalent} if and only if there exists a formal automorphism $\Phi$ such that the following two conditions hold:
   \begin{equation}
   \label{eq:equiv mc}
   \begin{split}
   \Dbar & \,=\, (\Phi^{-1})^{\otimes 2} \circ D \circ \Phi, \\
   \Lambdabar & \,=\, \Phi^{-1} \Lambda \Phi.
   \end{split}
   \end{equation}
The $H$-module-coalgebra $A$ is said to be \emph{rigid} if every deformation of $A$ is equivalent to the trivial deformation.

The following result is the module-coalgebra analogue of Theorem \ref{thm:def ma} and can be proved by similar arguments.

\begin{thm}
\label{thm:def mc}
Let $A$ be an $H$-module-coalgebra.  Then the following statements hold.
\begin{enumerate}
\item The infinitesimal $\theta_1$ of any deformation $\Theta$ of $A$ is a $2$-cocycle in the deformation complex $\Cmc^2(A)$ whose cohomology class is determined by the equivalence class of $\Theta$.
\item If $\Hmc^2(A) = 0$, then $A$ is rigid.
\end{enumerate}
\end{thm}


\section{Deformation bicomplex of comodule-algebras}
\label{sec:ca}

The purposes of this section are to construct the deformation bicomplex $\Cca^{**}(A) = \Hoch^*(A, H^{\otimes *} \otimes A)$ of an $H$-comodule-algebra $A$ and to list the corresponding deformation results.  Deformations are taken with respect to both the $H$-comodule structure and the algebra structure on $A$.

\subsection{Comodule-algebra}
\label{subsec:ca}

For an associative algebra $A = (A, \mu_A)$, an \emph{$H$-comodule-algebra} structure on $A$ consists of an $H$-comodule structure $\rho \in \Hom(A,H\otimes A)$ on $A$ such that the map $\mu_A \colon A \otimes A \to A$ becomes an $H$-comodule morphism, i.e.,
   \[
   \rho \circ \mu_A \,=\, (\mu_H \otimes \mu_A) \circ (\Id_H \otimes \tau_{(A,H)} \otimes \Id_A) \circ \rho^{\otimes 2}.
   \]
Here and in what follows, given two modules $X$ and $Y$,
   \[
   \tau_{(X,Y)} \colon X \otimes Y \cong Y \otimes X
   \]
denotes the twist isomorphism.  For the rest of this section, $A$ will denote an $H$-comodule-algebra with structure map $\rho$.

In general, for $n \geq 2$, the module $A^{\otimes n}$ becomes an $H$-comodule whose structure map
   \[
   \rho^n \colon A^{\otimes n} \to H \otimes A^{\otimes n}
   \]
is defined by the commutative diagram:
   \[
   \begin{CD}
   A^{\otimes n} @>{\rho^n}>> H \otimes A^{\otimes n} \\
   @V{\rho^{\otimes n}}VV     @AA{\mu_H^n \otimes \Id_{A^{\otimes n}}}A \\
   (H \otimes A)^{\otimes n} @>>{\mathrm{shuffle}}> H^{\otimes n} \otimes A^{\otimes n}.
   \end{CD}
   \]
Here
   \[
   \mu_H^n \colon H^{\otimes n} \to H
   \]
is the multiplication map defined by
   \begin{equation}
   \label{eq:mun}
   \mu_H^n(x_{1,n}) = x_1 \cdots x_n
   \end{equation}
for $x_i \in H$.

\subsection{$H^{\otimes q}$ as an $H$-bimodule}
\label{subsec:Hq}

For $q \geq 2$, the module $H^{\otimes q}$ is an $H$-bimodule whose left and right $H$-actions are defined by
   \[
   \begin{split}
   \mu_{l,H}^q(x, y_1 \otimes \cdots \otimes y_q) & \,=\, \sideset{}{_{(x)}}\sum x_{(1)}y_1 \otimes \cdots \otimes x_{(q)}y_q, \\
   \mu_{r,H}^q(y_1 \otimes \cdots \otimes y_q, x) & \,=\, \sideset{}{_{(x)}}\sum y_1x_{(1)} \otimes \cdots \otimes y_qx_{(q)}
   \end{split}
   \]
for $x, y_1, \ldots , y_q \in H$.

\subsection{$H^{\otimes q} \otimes A$ as an $A$-bimodule}
\label{subsec:HqA}

For each $q \geq 1$, the module $H^{\otimes q} \otimes A$ becomes an $A$-bimodule via the structure maps
   \[
   \begin{split}
   (\mu_{l,H}^q \otimes \mu_A) \circ (\Id_H \otimes \tau_{(A,H^{\otimes q})} \otimes A) \circ (\rho \otimes \Id_{H^{\otimes q} \otimes A}) & \colon A \otimes (H^{\otimes q} \otimes A) \to H^{\otimes q} \otimes A, \\
   (\mu_{r,H}^q \otimes \mu_A) \circ (\Id_{H^{\otimes q}} \otimes \tau_{(A,H)} \otimes \Id_A) \circ (\Id_{H^{\otimes q} \otimes A} \otimes \rho) & \colon  (H^{\otimes q} \otimes A) \otimes A \to  H^{\otimes q} \otimes A.
   \end{split}
   \]
In particular, we can consider the Hochschild cochain complex $\Hoch^*(A, H^{\otimes q} \otimes A)$.

\subsection{Comodule-algebra deformation bicomplex}
\label{subsec:bicomplex ca}

The \emph{deformation bicomplex} of the $H$-comodule-algebra $A$ is the first quadrant, cohomological bicomplex
   \[
   \Cca^{**}(A) \,=\, \left\lbrace \Cca^{p,\,q}(A), d^{p,\,q}, (-1)^{p+1}b^{p,\,q}\right\rbrace
   \]
with
   \[
   \Cca^{p,\,q}(A) \,=\,
   \begin{cases}
   0 & \text{ if } p = 0, \\
   \Der(A) & \text{ if } (p,q) = (1,0), \\
   \Hom(A^{\otimes p}, H^{\otimes q} \otimes A) & \text{ otherwise}.
   \end{cases}
   \]
In the $q$th row for $q \geq 0$, define the horizontal differential as the Hochschild coboundary,
   \[
   d^{*,\,q} \,\buildrel \text{def} \over=\, \delta_h^* \colon \Hoch^*(A, H^{\otimes q} \otimes A) \to \Hoch^{*+1}(A, H^{\otimes q} \otimes A).
   \]
For $p \geq 1$, the vertical differential
   \[
   (-1)^{p+1}b^{p,\,q} \colon \Hom(A^{\otimes p}, H^{\otimes q} \otimes A) \to \Hom(A^{\otimes p}, H^{\otimes (q+1)} \otimes A)
   \]
is defined by the alternating sum
   \[
   b^{p,\,q} \,=\, \sum_{i=0}^{q+1} \, (-1)^i b^{p,\,q} \lbrack i \rbrack,
   \]
where
   \[
   b^{p,\,q} \lbrack i \rbrack(\varphi) \,=\,
   \begin{cases}
   \left(\Id_{H^{\otimes q}} \otimes \rho\right) \circ \varphi & \text{ if } i = 0, \\
   \left(\Id_{H^{\otimes (q-i)}} \otimes \Delta_H \otimes \Id_{H^{\otimes (i-1)} \otimes A}\right) \circ \varphi & \text{ if } 1 \leq i \leq q, \\
   \left(\Id_H \otimes \varphi\right) \circ \rho^p & \text{ if } i = q+1.
   \end{cases}
   \]
An analogue of Theorem \ref{thm:bicomplex ma} holds in the comodule-algebra context.  The total complex of the deformation bicomplex $\Cca^{**}(A)$ is denoted by $\Cca^*(A)$ and is called the \emph{deformation complex} of $A$.  The $n$th cohomology module of $\Cca^*(A)$ is denoted by $\Hca^n(A)$.

Note that the cochain complex $\cF^*_{ca}(A)$ constructed in \cite{yau} is the total complex of the sub-bicomplex of $\Cca^{**}(A)$ consisting of the first column $\Cca^{1,\,*}(A)$ and the first row $\Cca^{*,\,1}(A)$.

\subsection{Cup-product}
\label{subsec:cup ca}

The module $H^{\otimes q} \otimes A$ is an associative algebra whose multiplication map is defined as
   \[
   \mu_{H^{\otimes q} \otimes A} \,\buildrel \text{def} \over =\,
   \left(\mu_{H^{\otimes q}} \otimes \mu_A \right) \circ
   \left(\Id_{H^{\otimes q}} \otimes \tau_{(A, H^{\otimes q})} \otimes \Id_A\right),
   \]
where
   \[
   \mu_{H^{\otimes q}}(x_1 \otimes \cdots \otimes x_q, y_1 \otimes \cdots \otimes y_q) \,=\,
   x_1y_1 \otimes \cdots \otimes x_qy_q
   \]
for $x_i, y_j \in H$.  The three conditions in \eqref{eq:3conditions} hold in this case.  It follows that each row $\Cca^{*,\,q}(A) = \Hoch^*(A, H^{\otimes q} \otimes A)$ admits a Hochschild $\cup$-product as in \eqref{eq:cup Hoch}.  Explicitly, given $\varphi \in \Cca^{r,\,q}(A)$ and $\psi \in \Cca^{s,\,q}(A)$, their $\cup$-product is given by
   \[
   \varphi \cup \psi \,=\, \mu_{H^{\otimes q} \otimes A} \circ (\varphi \otimes \psi).
   \]

\subsection{Comodule-algebra deformation}
\label{subsec:deformation ca}

A \emph{deformation} of $A$ is a power series $\Theta = \sum_{n \geq 0}\theta_n t^n$ with $\theta_0 = (\rho, \mu_A)$ and each $\theta_n = (\rho_n, \pi_n) \in \Cca^2(A)$, satisfying the following three conditions.  Writing $R = \sum_{n \geq 0} \rho_n t^n$ $(\rho_0 = \rho)$ and $\Pi = \sum_{n \geq 0} \pi_n t^n$ $(\pi_0 = \mu_A)$, the three conditions are:
   \begin{equation}
   \label{eq:def ca}
   \begin{split}
   (\Id_H \otimes R) \circ R & \,=\, (\Delta_H \otimes \Id_A) \circ R, \\
   R \circ \Pi & \,=\, (\mu_H \otimes \Pi) \circ (\Id_H \otimes \tau_{(A, H)} \otimes \Id_A) \circ R^{\otimes 2}, \\
   \Pi(\Pi(a,b),c) & \,=\, \Pi(a, \Pi(b,c))
   \end{split}
   \end{equation}
for $a, b, c \in A$.  Such a deformation is also denoted by $\Theta = (R, \Pi)$.  The linear coefficient $\theta_1 \in \Cca^2(A)$ is called the \emph{infinitesimal}.  The \emph{trivial deformation} is the deformation $\Theta = \theta_0 = (\rho, \mu_A)$.

A \emph{formal automorphism} of $A$ is a power series $\Phi = \sum_{n \geq 0} \phi_n t^n$ with $\phi_0 = \Id_A$ and each $\phi_n \in \End(A)$ such that the first non-zero $\phi_n$ $(n \geq 1)$ is a derivation on $A$.

Two deformations $\Theta = (R, \Pi)$ and $\Thetabar = (\Rbar, \Pibar)$ are said to be \emph{equivalent} if and only if there exists a formal automorphism $\Phi$ such that the following two conditions are satisfied:
   \begin{equation}
   \label{eq:equiv ca}
   \begin{split}
   \Rbar & \,=\, (\Id_H \otimes \Phi^{-1}) \circ R \circ \Phi, \\
   \Pibar & \,=\, \Phi^{-1} \circ \Pi \circ \Phi^{\otimes 2}.
   \end{split}
   \end{equation}
The $H$-comodule-algebra $A$ is said to be \emph{rigid} if every deformation of $A$ is equivalent to the trivial deformation.

The following result is the comodule-algebra analogue of Theorem \ref{thm:def ma}.

\begin{thm}
\label{thm:def ca}
Let $A$ be an $H$-comodule-algebra.  Then the following statements hold.
\begin{enumerate}
\item The infinitesimal $\theta_1$ of any deformation $\Theta$ of $A$ is a $2$-cocycle in the deformation complex $\Cca^2(A)$ whose cohomology class is determined by the equivalence class of $\Theta$.
\item If $\Hca^2(A) = 0$, then $A$ is rigid.
\end{enumerate}
\end{thm}


\section{Deformation bicomplex of comodule-coalgebras}
\label{sec:cc}

The purposes of this section are to construct the deformation bicomplex $\Ccc^{**}(A) = \Hom(A, H^{\otimes *} \otimes A^{\otimes *})$ of an $H$-comodule-coalgebra $A$ and to list the corresponding deformation results.  Deformations are taken with respect to both the $H$-comodule structure and the coalgebra structure on $A$.

\subsection{Comodule-coalgebra}
\label{subsec:cc}

For this section, let $A = (A, \Delta_A)$ be a coassociative coalgebra.  An \emph{$H$-comodule-coalgebra} structure on $A$ consists of an $H$-comodule structure $\rho \in \Hom(A, H \otimes A)$ on $A$ such that the comultiplication map $\Delta_A$ is an $H$-comodule morphism, i.e.,
   \[
   (\Id_H \otimes \Delta_A) \circ \rho
   \,=\,
   (\mu_H \otimes \Id_{A^{\otimes 2}}) \circ (\Id_H \otimes \tau_{(A,H)} \otimes \Id_A) \circ \rho^{\otimes 2} \circ \Delta_A.
   \]
For the rest of this section, $A$ will denote an $H$-comodule-coalgebra with structure map $\rho$.

\subsection{Comodule-coalgebra deformation bicomplex}
\label{subsec:bicomplex cc}

The \emph{deformation bicomplex} of $A$ is the first quadrant, cohomological bicomplex
   \[
   \Ccc^{**}(A) \,=\, \left\lbrace\Ccc^{p,\,q}(A), d^{p,\,q}, (-1)^{p+1}b^{p,\,q}\right\rbrace
   \]
with
   \[
   \Ccc^{p,\,q}(A) \,=\,
   \begin{cases}
   0 & \text{ if } p = 0, \\
   \Coder(A) & \text{ if } (p,q) = (1,0), \\
   \Hom(A, H^{\otimes q} \otimes A^{\otimes p}) & \text{ otherwise}.
   \end{cases}
   \]
For $p \geq 1$, the vertical differential
   \[
   (-1)^{p+1}b^{p,\,q} \colon \Hom(A, H^{\otimes q} \otimes A^{\otimes p}) \to \Hom(A, H^{\otimes (q+1)} \otimes A^{\otimes p})
   \]
is defined by the alternating sum
   \[
   b^{p,\,q} \,=\, \sum_{i=0}^{q+1} \, (-1)^i b^{p,\,q} \lbrack i \rbrack,
   \]
where
   \[
   b^{p,\,q}\lbrack i \rbrack(\varphi) \,=\,
   \begin{cases}
   \left(\Id_{H^{\otimes q}} \otimes \rho^p\right) \circ \varphi & \text{ if } i = 0, \\
   \left(\Id_{H^{\otimes (q-i)}} \otimes \Delta_H \otimes \Id_{H^{\otimes (i-1)} \otimes A^{\otimes p}} \right) \circ \varphi & \text{ if } 1 \leq i \leq q, \\
   \left(\Id_H \otimes \varphi\right) \circ \rho & \text{ if } i = q+1
   \end{cases}
   \]
for $\varphi \in \Hom(A, H^{\otimes q} \otimes A^{\otimes p})$.  Here $\rho^p$ is the $H$-comodule structure map on $A^{\otimes p}$ defined in section \ref{subsec:ca}.

In the $0$th row $\Ccc^{*,\,0}(A)$, define the horizontal differential as the Hochschild coalgebra coboundary,
   \[
   d^{p,\,0} \,\buildrel \text{def} \over=\, \delta^p_c \colon \Hoch^p_c(A,A) \to \Hoch^{p+1}_c(A,A).
   \]
In the $q$th row $\Ccc^{*,\,q}(A)$ for $q \geq 1$, define the horizontal differential as the alternating sum
   \[
   d^{p,\,q} \,=\, \sum_{i=0}^{p+1} \, (-1)^i d^{p,\,q}\lbrack i \rbrack,
   \]
where
   \[
   d^{p,\,q} \lbrack i \rbrack(\varphi)
   \,=\,
   \begin{cases}
   \left(\mu^q_{l,H} \otimes \Id_{A^{\otimes (p+1)}}\right) \circ \left(\Id_H \otimes \tau_{(A,H^{\otimes q})} \otimes \Id_{A^{\otimes p}}\right) \circ (\rho \otimes \varphi) \circ \Delta_A & \text{ if } i =0, \\
   \left(\Id_{H^{\otimes q} \otimes A^{\otimes (i-1)}} \otimes \Delta_A \otimes \Id_{A^{\otimes (p-i)}}\right) \circ \varphi & \text{ if } 1 \leq i \leq p, \\
   \left(\mu^q_{r,H} \otimes \Id_{A^{\otimes (p+1)}}\right) \circ \left(\Id_{H^{\otimes q}} \otimes \tau_{(A^{\otimes p}, H)} \otimes \Id_A\right) \circ (\varphi \otimes \rho) \circ \Delta_A & \text{ if } i = p+1
   \end{cases}
   \]
for $\varphi \in \Hom(A, H^{\otimes q} \otimes A^{\otimes p})$.  Here $\mu^q_{*,H}$ $(* = l, r)$ are the $H$-bimodule structure maps on $H^{\otimes q}$ defined in section \ref{subsec:Hq}.

The comodule-coalgebra analogue of Theorem \ref{thm:bicomplex ma} holds.  The total complex of the deformation bicomplex $\Ccc^{**}(A)$ is denoted by $\Ccc^*(A)$ and is called the \emph{deformation complex} of $A$.  The $n$th cohomology module of $\Ccc^*(A)$ is denoted by $\Hcc^n(A)$.

Note that the cochain complex $\cF^*_{cc}(A)$ constructed in \cite{yau} is the total complex of the sub-bicomplex of $\Ccc^{**}(A)$ consisting of the first column $\Ccc^{1,\,*}(A)$ and the first row $\Ccc^{*,\,1}(A)$.

\subsection{Cup-product}
\label{subsec:cup cc}

For each $q \geq 0$, the $q$th row $(\Ccc^{*,\,q}(A), d^{*,\,q})$ is a DGA whose product is defined as
   \[
   \varphi \cup \psi \,\buildrel \text{def} \over=\,
   \left(\mu_{H^{\otimes q}} \otimes \Id_{A^{\otimes (r+s)}}\right) \circ \left(\Id_{H^{\otimes q}} \otimes \tau_{(A^{\otimes r}, H^{\otimes q})} \otimes \Id_{A^{\otimes s}}\right) \circ (\varphi \otimes \psi) \circ \Delta_A
   \]
for $\varphi \in \Ccc^{r,\,q}(A)$ and $\psi \in \Ccc^{s,\,q}(A)$.

\subsection{Comodule-coalgebra deformation}
\label{subsec:deformation cc}

A \emph{deformation} of $A$ is a power series $\Theta = \sum_{n \geq 0} \theta_n t^n$ with $\theta_0 = (\rho, \Delta_A)$ and each $\theta_n = (\rho_n, \Delta_n) \in \Ccc^2(A)$, satisfying the following three conditions.  Writing $R = \sum_{n \geq 0} \rho_n t^n$ $(\rho_0 = \rho)$ and $D = \sum_{n \geq 0} \Delta_n t^n$ $(\Delta_0 = \Delta_A)$, the three conditions are:
   \begin{equation}
   \label{eq:def cc}
   \begin{split}
   (\Id_H \otimes R) \circ R & \,=\, (\Delta_H \otimes \Id_A) \circ R, \\
   (\Id_H \otimes D) \circ R & \,=\, \left(\mu_H \otimes \Id_{A^{\otimes 2}}\right) \circ \left(\Id_H \otimes \tau_{(A,H)} \otimes \Id_A\right) \circ R^{\otimes 2} \circ D, \\
   (D \otimes \Id_A) \circ D & \,=\, (\Id_A \otimes D) \circ D.
   \end{split}
   \end{equation}
Such a deformation is also denoted by $\Theta = (R, D)$.  The linear coefficient $\theta_1 \in \Ccc^2(A)$ is called the \emph{infinitesimal}.  The \emph{trivial deformation} is the deformation $\Theta = \theta_0 = (\rho, \Delta_A)$.

A \emph{formal automorphism} of $A$ is a power series $\Phi = \sum_{n \geq 0} \phi_n t^n$ with $\phi_0 = \Id_A$ and each $\phi_n \in \End(A)$ such that the first non-zero $\phi_n$ $(n \geq 1)$ is a coderivation on $A$.

Two deformations $\Theta = (R, D)$ and $\Thetabar = (\Rbar, \Dbar)$ are said to be \emph{equivalent} if and only if there exists a formal automorphism $\Phi$ such that the following two conditions are satisfied:
   \begin{equation}
   \label{eq:equiv cc}
   \begin{split}
   \Rbar & \,=\, (\Id_H \otimes \Phi^{-1}) \circ R \circ \Phi, \\
   \Dbar & \,=\, (\Phi^{-1})^{\otimes 2} \circ D \circ \Phi.
   \end{split}
   \end{equation}
The $H$-comodule-coalgebra $A$ is said to be \emph{rigid} if every deformation of $A$ is equivalent to the trivial deformation.

The following result is the comodule-coalgebra analogue of Theorem \ref{thm:def ma}.

\pagebreak

\begin{thm}
\label{thm:def cc}
Let $A$ be an $H$-comodule-coalgebra.  Then the following statements hold.
\begin{enumerate}
\item The infinitesimal $\theta_1$ of any deformation $\Theta$ of $A$ is a $2$-cocycle in the deformation complex $\Ccc^2(A)$ whose cohomology class is determined by the equivalence class of $\Theta$.
\item If $\Hcc^2(A) = 0$, then $A$ is rigid.
\end{enumerate}
\end{thm}


\section{Deformation tricomplex of module-bialgebras}
\label{sec:mb}

In this section, we construct the deformation tricomplex $\Cmb^{***}(A) = \Hoch^*(H, \Hom(A^{\otimes *}, A^{\otimes *}))$ of an $H$-module-bialgebra $A$ and list the corresponding deformation results.  Both the $H$-module structure and the bialgebra structure on $A$ are deformed.  This deformation tricomplex contains the deformation bicomplexes $\Cma^{**}(A)$ and $\Cmc^{**}(A)$ for module-(co)algebras and the Gerstenhaber-Schack bicomplex $\CGS^{**}(A)$ \cite{gs1,gs2} for a bialgebra.  Unlike the (co)module-(co)algebra deformation bicomplexes, each cochain complex obtained from $\Cmb^{***}(A)$ by fixing two of the three dimensions is a DGA. In fact, each such cochain complex is either a Hochschild (coalgebra) cochain complex or is isomorphic to one. The relevant products can be identified with the Hochschild (coalgebra) $\cup$-products.

For the rest of this paper, let $A = (A, \mu_A, \Delta_A)$ be a bialgebra with associative multiplication $\mu_A$ and coassociative comultiplication $\Delta_A$.

\subsection{Biderivation}
\label{subsec:bider}

A \emph{biderivation on $A$} is a linear self-map $\varphi \in \End(A)$ that is both a derivation and a coderivation on $A$.  The set of biderivations on $A$ is denoted by $\Bider(A)$, and it is considered as a submodule of $\End(A)$.

\subsection{Module-bialgebra}
\label{subsec:module-bialgebra}

An \emph{$H$-module-bialgebra} structure on $A$ is an $H$-module structure $\lambda \in \Hom(H, \End(A))$ on $A$ that makes $A$ into an $H$-module-algebra and an $H$-module-coalgebra simultaneously.  For the rest of this section, let $A$ be an $H$-module-bialgebra with structure map $\lambda$.

For example, let $G_1$ and $G_2$ be two groups.  Then any group homomorphism
   \[
   \phi \colon G_1 \to \Aut(G_2)
   \]
gives rise to a $K \lbrack G_1 \rbrack$-module-bialgebra structure on the group bialgebra $K \lbrack G_2 \rbrack$ via the action
   \begin{equation}
   \label{eq:mb example}
   \lambda(x)(y) \,=\, \phi(x)(y)
   \end{equation}
for $x \in G_1$ and $y \in G_2$ \cite[Example 3.6]{abe}.  Similarly, suppose that $L_1$ and $L_2$ are Lie algebras and that
   \[
   \phi \colon L_1 \to \Der(L_2)
   \]
is a Lie algebra morphism.  Then the same formula \eqref{eq:mb example} gives rise to a $U(L_1)$-module-bialgebra structure on the enveloping bialgebra $U(L_2)$ \cite[Exampl 3.7]{abe}.

\subsection{$H$-bimodule structure on $\Hom(A^{\otimes p}, A^{\otimes q})$}
\label{subsec:Hom(Ap,Aq)}

When $A$ is an $H$-module-bialgebra, there is an $H$-bimodule structure on the module $\Hom(A^{\otimes p}, A^{\otimes q})$ for $p, q \geq 1$.  The left and right $H$-actions are given as follows:
   \[
   \begin{split}
   (x \varphi)(a_{1,p})
   & \,=\, \sum \lambda(x_{(1)})(\varphi(a_{1,p})^1) \otimes \cdots \otimes \lambda(x_{(q)})(\varphi(a_{1,p})^q), \\
   (\varphi x)(a_{1,p})
   & \,=\, \sum \varphi\left(\lambda(x_{(1)})(a_1) \otimes \cdots \otimes \lambda(x_{(p)})(a_p)\right).
   \end{split}
   \]
Here $x \in H$, $\varphi \in \Hom(A^{\otimes p}, A^{\otimes q})$, $a_i \in A$, and $\varphi(a_{1,p}) = \sum \varphi(a_{1,p})^1 \otimes \cdots \otimes \varphi(a_{1,p})^q$.  This generalizes the constructions in sections \ref{subsec:H-bimod} and \ref{subsec:Hom(A,An)}.  In particular, we can consider the Hochschild cochain complex $\Hoch^*(H, \Hom(A^{\otimes p}, A^{\otimes q}))$.

\subsection{Module-bialgebra deformation tricomplex}
\label{subsec:def tri mb}

The deformation tricomplex of the $H$-module-bialgebra $A$ is the first octant, cohomological tricomplex
   \[
   \Cmb^{***}(A) \,=\,
   \left\lbrace \Cmb^{p,\,q,\,r}(A),\, (-1)^{q+1}\dI^{p,\,q,\,r},\, (-1)^{r+1}\dII^{p,\,q,\,r},\, (-1)^{p+1}\dIII^{p,\,q,\,r} \right\rbrace
   \]
with
   \[
   \Cmb^{p,\,q,\,r}(A) \,=\,
   \begin{cases}
   0 & \text{ if $p = 0$ or $q = 0$}, \\
   \Bider(A) & \text{ if $(p, q, r) = (1, 1, 0)$}, \\
   \Hom(H^{\otimes r}, \Hom(A^{\otimes p}, A^{\otimes q})) & \text{ otherwise}.
   \end{cases}
   \]
The differential
   \[
   (-1)^{q+1} \dI^{p,\,q,\,r} \colon \Cmb^{p,\,q,\,r}(A) \to \Cmb^{p+1,\, q,\,r}(A)
   \]
is defined using the alternating sum
   \[
   \dI^{p,\,q,\,r} \,=\, \sum_{i=0}^{p+1} \, (-1)^i \dI^{p,\,q,\,r} \lbrack i \rbrack,
   \]
where
   \[
   \left(\dI^{p,\,q,\,r}\lbrack i \rbrack \varphi\right)(x_{1,r})(a_{1,p+1}) \,=\,
   \begin{cases}
   \sum\mu_{l,A}^q\left(\lambda(x_{1(1)} \cdots x_{r(1)})(a_1), \varphi(x_{1(2)} \otimes \cdots \otimes x_{r(2)})(a_{2,p+1})\right) & \text{ if $i = 0$}, \\
   \varphi(x_{1,r})(a_{1,i-1} \otimes (a_i a_{i+1}) \otimes a_{i+2, p+1}) & \text{ if $1 \leq i \leq p$}, \\
   \sum\mu_{r,A}^q\left(\varphi(x_{1(1)} \otimes \cdots \otimes x_{r(1)})(a_{1,p}), \lambda(x_{1(2)} \cdots x_{r(2)})(a_{p+1})\right) & \text{ if $i = p+1$}
   \end{cases}
   \]
for $\varphi \in \Cmb^{p,\,q,\,r}(A)$, $x_k \in H$, and $a_l \in A$.   Here $\mu_{l,A}^q$ and $\mu_{r,A}^q$ are the left and right actions of $A$ on $A^{\otimes q}$, defined as in section \ref{subsec:Hq}.

The differential
   \[
   (-1)^{r+1} \dII^{p,\,q,\,r} \colon \Cmb^{p,\,q,\,r}(A) \to \Cmb^{p,\,q+1,\,r}(A)
   \]
is defined using the alternating sum
   \[
   \dII^{p,\,q,\,r} \,=\, \sum_{i=0}^{q+1} \, (-1)^i \dII^{p,\,q,\,r} \lbrack i \rbrack,
   \]
where
   \begin{multline*}
   \left(\dII^{p,\,q,\,r} \lbrack i \rbrack \varphi\right)(x_{1,r})(a_{1,p}) \\
   \,=\,
   \begin{cases}
   \sum \lambda(x_{1(1)} \cdots x_{r(1)})(a_{1(1)} \cdots a_{p(1)}) \otimes \varphi(x_{1(2)} \otimes \cdots \otimes x_{r(2)})(a_{1(2)} \otimes \cdots \otimes a_{p(2)}) & \text{ if $i = 0$}, \\
   \left(\Id_{A^{\otimes (i-1)}} \otimes \Delta_A \otimes \Id_{A^{\otimes (q-i)}}\right) \left(\varphi(x_{1,r})(a_{1,p})\right) & \text{ if $1 \leq i \leq q$}, \\
   \sum \varphi(x_{1(1)} \otimes \cdots \otimes x_{r(1)})(a_{1(1)} \otimes \cdots \otimes a_{p(1)}) \otimes \lambda(x_{1(2)} \cdots x_{r(2)})(a_{1(2)} \cdots a_{p(2)}) & \text{ if $i = q+1$}.
   \end{cases}
   \end{multline*}

Finally, the differential
   \[
   (-1)^{p+1} \dIII^{p,\,q,\,r} \colon \Cmb^{p,\,q,\,r}(A) \to \Cmb^{p,\,q,\,r+1}(A)
   \]
is defined using the Hochschild coboundary
   \[
   \dIII^{p,\,q,\,r}
   \buildrel \text{def} \over =\,
   \delta^r_h \colon
   \Hoch^r(H, \Hom(A^{\otimes p}, A^{\otimes q})) \to
   \Hoch^{r+1}(H, \Hom(A^{\otimes p}, A^{\otimes q})).
   \]
In particular, $\dIII^{p,\,q,\,r}$ is the alternating sum
   \[
   \dIII^{p,\,q,\,r} \,=\, \sum_{i=0}^{r+1} \, (-1)^i \dIII^{p,\,q,\,r} \lbrack i \rbrack,
   \]
where
   \[
   \dIII^{p,\,q,\,r} \lbrack i \rbrack \,=\, \delta^r_h \lbrack i \rbrack
   \]
as in \eqref{eq:Hoch diff}.

The following observations, which are the module-bialgebra analogue of Theorem \ref{thm:bicomplex ma}, ensure that $\Cmb^{***}(A)$ is indeed a tricomplex.

\begin{thm}
\label{thm:bicomplex mb}
Let $A$ be an $H$-module-bialgebra.  Then:
\begin{enumerate}
\item The following statements hold for all possible values of $p, q, r$ and $k < l$:
   \[
   \begin{split}
   \dI^{p+1,\,q,\,r} \lbrack l \rbrack \circ \dI^{p,\,q,\,r} \lbrack k \rbrack
   & \,=\, \dI^{p+1,\,q,\,r} \lbrack k \rbrack \circ \dI^{p,\,q,\,r} \lbrack l - 1 \rbrack, \\
   \dII^{p,\,q+1,\,r} \lbrack l \rbrack \circ \dII^{p,\,q,\,r} \lbrack k \rbrack
   & \,=\, \dII^{p,\,q+1,\,r} \lbrack k \rbrack \circ \dII^{p,\,q,\,r} \lbrack l - 1 \rbrack, \\
   \dIII^{p,\,q,\,r+1} \lbrack l \rbrack \circ \dIII^{p,\,q,\,r} \lbrack k \rbrack
   & \,=\, \dIII^{p,\,q,\,r+1} \lbrack k \rbrack \circ \dIII^{p,\,q,\,r} \lbrack l - 1 \rbrack.
   \end{split}
   \]
\item The following statements hold for all possible values of $p, q, r, i$, and $j$:
   \[
   \begin{split}
   \dII^{p+1,\,q,\,r} \lbrack j \rbrack \circ \dI^{p,\,q,\,r} \lbrack i \rbrack
   & \,=\, \dI^{p,\,q+1,\,r} \lbrack i \rbrack \circ \dII^{p,\,q,\,r} \lbrack j \rbrack, \\
   \dII^{p,\,q,\,r+1} \lbrack j \rbrack \circ \dIII^{p,\,q,\,r} \lbrack i \rbrack
   & \,=\, \dIII^{p,\,q+1,\,r} \lbrack i \rbrack \circ \dII^{p,\,q,\,r} \lbrack j \rbrack, \\
   \dIII^{p+1,\,q,\,r} \lbrack j \rbrack \circ \dI^{p,\,q,\,r} \lbrack i \rbrack
   & \,=\, \dI^{p,\,q,\,r+1}\lbrack i \rbrack \circ \dIII^{p,\,q,\,r} \lbrack j \rbrack.
   \end{split}
   \]
\end{enumerate}
In particular, every two-dimensional plane in $\Cmb^{***}(A)$ is a bicomplex.
\end{thm}

\subsection{Boundary planes}
\label{subsec:boundary plane}

Observe that the boundary planes of the deformation tricomplex $\Cmb^{***}(A)$ are either known or have been discussed in previous sections.  In fact:
\begin{enumerate}
\item The $p = 1$ plane $\Cmb^{1,\,*,\,*}(A)$ coincides with the deformation bicomplex $\Cmc^{**}(A)$ in which $A$ is regarded as an $H$-module-coalgebra (see section \ref{sec:mc}).
\item The $q = 1$ plane $\Cmb^{*,\,1,\,*}(A)$ coincides with the deformation bicomplex $\Cma^{**}(A)$ in which $A$ is regarded as an $H$-module-algebra (see section \ref{sec:bicomplex}).
\item The $r = 0$ plane $\Cmb^{*,\,*,\,0}(A)$ coincides with the Gerstenhaber-Schack deformation bicomplex $\CGS^{**}(A)$ (denoted by $\widehat{C}_b^{\bullet,\, \bullet}(A,A)$ in \cite{gs1,gs2}) in which $A$ is regarded as only a bialgebra.
\end{enumerate}
As in (co)module-(co)algebra cases, the only exception to the above remarks is the entry $\Cmb^{1,\,1,\,0}(A) = \Bider(A)$.

\subsection{Cup-products}
\label{subsec:cup mb}

Each cochain complex in the deformation tricomplex $\Cmb^{***}(A)$ has an associative $\cup$-product that makes it into a DGA.

In the direction of $\dI$, fix $q \geq 1$ and $r \geq 0$.  Then the cochain complex $\Cmb^{*,\,q,\,r}(A)$ is a DGA whose product is defined as
   \begin{equation}
   \label{eq:cup Cqr}
   (\varphi \cup \psi)(x_{1,r})(a_{1,p_1+p_2}) \,\buildrel \text{def} \over =\,
   \mu_{A^{\otimes q}}\left(\varphi(x_{1(1)} \otimes \cdots \otimes x_{r(1)})(a_{1,p_1}),\,
   \psi(x_{1(2)} \otimes \cdots \otimes x_{r(2)})(a_{p_1+1, p_1+p_2})\right)
   \end{equation}
for $\varphi \in \Cmb^{p_1,\,q,\,r}(A)$ and $\psi \in \Cmb^{p_2,\,q,\,r}$.  This generalizes the $\cup$-product in $\Cma^{*,\,r}(A)$ discussed in section \ref{sec:cup}.

In fact, it can be identified with a Hochschild $\cup$-product.  There is a canonical isomorphism of cochain complexes,
   \begin{equation}
   \label{eq:Cqr}
   \left(\Cmb^{*,\,q,\,r}(A), (-1)^{q+1} \dI^{*,\,q,\,r}\right)
   \,\cong\,
   (-1)^{q+1} \left(\Hoch^*(A, \Hom(H^{\otimes r}, A^{\otimes q})), \delta_h\right)
   \end{equation}
given by the $\Hom-\otimes$ adjunction (twice).  The left and right actions of $A$ on $\Hom(H^{\otimes r}, A^{\otimes q})$ are given by
   \[
   \begin{split}
   (af)(x_{1,r}) & \,=\, \sum \mu^q_{l,A}\left(\lambda(x_{1(1)} \cdots x_{r(1)})(a), f(x_{1(2)} \otimes \cdots \otimes x_{r(2)})\right), \\
   (fa)(x_{1,r}) & \,=\, \sum \mu^q_{r,A}\left(f(x_{1(1)} \otimes \cdots \otimes x_{r(1)}), \lambda(x_{1(2)} \cdots x_{r(2)})(a)\right)
   \end{split}
   \]
for $a \in A$, $f \in \Hom(H^{\otimes r}, A^{\otimes q})$, and $x_i \in H$.  The module $\Hom(H^{\otimes r}, A^{\otimes q})$ is an associative algebra via the product
   \[
   (f \cdot g)(x_{1,r})
   \,=\,
   \mu_{A^{\otimes q}}\left(f(x_{1(1)} \otimes \cdots \otimes x_{r(1)}), g(x_{1(2)} \otimes \cdots \otimes x_{r(2)})\right)
   \]
such that the conditions \eqref{eq:3conditions} are satisfied.  The resulting Hochschild $\cup$-product corresponds, via the isomorphism \eqref{eq:Cqr}, to the $\cup$-product defined in \eqref{eq:cup Cqr}.

In the direction of $\dII$, fix $p \geq 1$ and $r \geq 0$.  Then the cochain complex $\Cmb^{p,\,*,\,r}(A)$ is a DGA whose product is defined as
   \begin{multline}
   \label{eq:cup Cpr}
   (\varphi \cup \psi)(x_{1,r})(a_{1,p}) \,\buildrel \text{def} \over =\, \\
   \sum \varphi(x_{1(1)} \otimes \cdots \otimes x_{r(1)})(a_{1(1)} \otimes \cdots \otimes a_{p(1)}) \otimes \psi(x_{1(2)} \otimes \cdots \otimes x_{r(2)})(a_{1(2)} \otimes \cdots \otimes a_{p(2)})
   \end{multline}
for $\varphi \in \Cmb^{p,\,q_1,\,r}(A)$ and $\psi \in \Cmb^{p,\,q_2,\,r}(A)$.  This generalizes the $\cup$-product in $\Cmc^{*,\,r}(A)$ discussed in section \ref{subsec:cup mc}.

Moreover, there is a canonical isomorphism,
   \begin{equation}
   \label{eq:Cpr}
   \left(\Cmb^{p,\,*,\,r}(A), (-1)^{r+1}\dII^{p,\,*,\,r}\right)
   \,\cong\,
   (-1)^{r+1} \left(\Hoch^*_c(H^{\otimes r} \otimes A^{\otimes p}, A), \delta_c\right)
   \end{equation}
of cochain complexes given by the $\Hom-\otimes$ adjunction.  The left and right $A$-coactions on $H^{\otimes r} \otimes A^{\otimes p}$ are given by
   \[
   \begin{split}
   x_{1,r} \otimes a_{1,p} & \,\mapsto\, \lambda(x_{1(1)} \cdots x_{r(1)})(a_{1(1)} \cdots a_{p(1)}) \otimes \left(x_{1(2)} \otimes \cdots \otimes x_{r(2)} \otimes a_{1(2)} \otimes \cdots \otimes a_{p(2)}\right), \\
   x_{1,r} \otimes a_{1,p} & \,\mapsto\, \left(x_{1(1)} \otimes \cdots \otimes x_{r(1)} \otimes a_{1(1)} \otimes \cdots \otimes a_{p(1)}\right) \otimes \lambda(x_{1(2)} \cdots x_{r(2)})(a_{1(2)} \cdots a_{p(2)}).
   \end{split}
   \]
The module $H^{\otimes r} \otimes A^{\otimes p}$ is a coassociative coalgebra via the coproduct
   \[
   \Delta(x_{1,r} \otimes a_{1,p})
   \,=\,
   \sum \left(x_{1(1)} \otimes \cdots \otimes x_{r(1)} \otimes a_{1(1)} \otimes \cdots \otimes a_{p(1)}\right) \otimes \left(x_{1(2)} \otimes \cdots \otimes x_{r(2)} \otimes a_{1(2)} \otimes \cdots \otimes a_{p(2)}\right)
   \]
such that the three conditions in \eqref{eq:3cond} are satisfied.  The resulting Hochschild coalgebra $\cup$-product corresponds, via the isomorphism \eqref{eq:Cpr}, to the $\cup$-product defined in \eqref{eq:cup Cpr}.

In the direction of $\dIII$, fix $p \geq 1$ and $q \geq 1$.  Note that the module $\Hom(A^{\otimes p}, A^{\otimes q})$ is an associative algebra via the product
   \[
   f \cdot g \,\buildrel \text{def} \over =\, f \circ \Delta_A^{p-1} \circ \mu_A^q \circ g,
   \]
where $\Delta_A^0 = \Id_A = \mu_A^1$.  In the particular case $p = q = 1$, this product is simply the composition product in $\Hom(A,A)$.  Regarding $\Hom(A^{\otimes p}, A^{\otimes q})$ as an $H$-bimodule, the three conditions in \eqref{eq:3conditions} hold.  It follows that the Hochschild cochain complex $\Cmb^{p,\,q,\,*}(A) = (-1)^{p+1}\Hoch^*(H,\Hom(A^{\otimes p}, A^{\otimes q}))$ admits a Hochschild $\cup$-product \eqref{eq:cup Hoch} that makes it into a DGA.

\subsection{Total complex}
\label{subsec:tot mb}

The total complex $\Cmb^*(A)$, called the \emph{deformation complex of $A$}, of the deformation tricomplex $\Cmb^{***}(A)$ is defined as usual but with a shift of degree:
   \[
   \Cmb^n(A) \,\buildrel \text{def} \over=\, \bigoplus_{p+q+r \,=\, n+1} \Cmb^{p,\,q,\,r}(A).
   \]
In particular, we have that
   \[
   \begin{split}
   \Cmb^1(A) & \,=\, \Bider(A), \\
   \Cmb^2(A) & \,=\, \Hom(H, \Hom(A,A)) \oplus \Hom(A^{\otimes 2}, A) \oplus \Hom(A, A^{\otimes 2}).
   \end{split}
   \]
The degree shift is introduced to ensure that the deformation results below have the same degree conventions as in the previous sections and as in \cite{gs1,gs2}.

The $n$th cohomology module of the deformation complex $\Cmb^*(A)$ is denoted by $\Hmb^n(A)$.

\subsection{Module-bialgebra deformation}
\label{subsec:def mb}

A \emph{deformation} of $A$ as an $H$-module-bialgebra is a power series $\Theta = \sum_{n \geq 0} \theta_n t^n$ with $\theta_0 = (\lambda, \mu_A, \Delta_A) \in \Cmb^2(A)$ and each $\theta_n = (\lambda_n, \pi_n, \Delta_n) \in \Cmb^2(A)$, satisfying the five conditions for module-algebra deformations \eqref{eq:def ma axioms} and module-coalgebra deformations \eqref{eq:def mc}.  Using the notations from earlier sections, such a deformation is also denoted by $\Theta = (\Lambda, \Pi, D)$.  The \emph{infinitesimal} of a deformation $\Theta$ is the linear coefficient $\theta_1$.

A \emph{formal automorphism} of $A$ is a power series $\sum_{n \geq 0} \phi_n t^n$ with $\phi_0 = \Id_A$ and each $\phi_n \in \End(A)$ such that the first non-zero $\phi_n$ $(n \geq 1)$ is a biderivation on $A$.

Two deformations $\Theta = (\Lambda, \Pi, D)$ and $\Thetabar = (\Lambdabar, \Pibar, \Dbar)$ are said to be \emph{equivalent} if and only if there exists a formal automorphism $\Phi$ such that the three conditions in \eqref{eq:equiv ma} and \eqref{eq:equiv mc} are satisfied.

The \emph{trivial deformation} is the deformation $\Theta = \theta_0 = (\lambda, \mu_A, \Delta_A)$.  The $H$-module-bialgebra $A$ is said to be \emph{rigid} if every deformation of $A$ is equivalent to the trivial deformation.

The following result is the module-bialgebra analogue of Theorem \ref{thm:def ma}.

\begin{thm}
\label{thm:def mb}
Let $A$ be an $H$-module-bialgebra.  Then the following statements hold.
\begin{enumerate}
\item The infinitesimal $\theta_1$ of any deformation $\Theta$ of $A$ is a $2$-cocycle in the deformation complex $\Cmb^2(A)$ whose cohomology class is determined by the equivalence class of $\Theta$.
\item If $\Hmb^2(A) = 0$, then $A$ is rigid.
\end{enumerate}
\end{thm}


\section{Deformation tricomplex of comodule-bialgebras}
\label{sec:cb}

The purposes of this section are to construct the deformation tricomplex $\Ccb^{***}(A) = \Hom(A^{\otimes *}, H^{\otimes *} \otimes A^{\otimes *})$ of an $H$-comodule-bialgebra $A$ and to list the corresponding deformation results.  Both the $H$-comodule structure and the bialgebra structure on $A$ are deformed.  This deformation tricomplex contains the deformation bicomplexes $\Cca^{**}(A)$ and $\Ccc^{**}(A)$ for comodule-(co)algebras and the Gerstenhaber-Schack bicomplex $\CGS^{**}(A)$ \cite{gs1,gs2} for a bialgebra.

\subsection{Comodule-bialgebra}
\label{subsec:comodule-bialgebra}

An \emph{$H$-comodule-bialgebra} structure on a bialgebra $A$ is an $H$-comodule structure $\rho \in \Hom(A, H \otimes A)$ that makes $A$ into an $H$-comodule-algebra and an $H$-comodule-coalgebra simultaneously.

For example, let $H$ be a commutative Hopf algebra with antipode $S$.  Then the map $\rho \in \Hom(H, H \otimes H)$ defined by
   \[
   \rho(x) \,=\, \sum \left(x_{(1)} Sx_{(3)}\right) \otimes x_{(2)}
   \]
for $x \in H$ gives $H$ an $H$-comodule-bialgebra structure \cite[Example 3.8]{abe}.

\subsection{$A$-bimodule structure on $H^{\otimes r} \otimes A^{\otimes q}$}
\label{subsec:HrAq}

For $r \geq 0$ and $q \geq 1$, there is an $A$-bimodule structure on the module $H^{\otimes r} \otimes A^{\otimes q}$ whose left and right $A$-action maps are defined as
   \[
   \begin{split}
   \left(\mu_{l,H}^r \otimes \mu_{l,A}^q\right) \circ \left(\Id_H \otimes \tau_{(A, H^{\otimes r})} \otimes \Id_{A^{\otimes q}}\right) \circ \left(\rho \otimes \Id_{H^{\otimes r} \otimes A^{\otimes q}}\right)
   & \colon A \otimes (H^{\otimes r} \otimes A^{\otimes q}) \to H^{\otimes r} \otimes A^{\otimes q}, \\
   \left(\mu_{r,H}^r \otimes \mu_{r,A}^q\right) \circ \left(\Id_{H^{\otimes r}} \otimes \tau_{(A^{\otimes q}, H)} \otimes \Id_A\right) \circ \left(\Id_{H^{\otimes r} \otimes A^{\otimes q}} \otimes \rho\right) & \colon (H^{\otimes r} \otimes A^{\otimes q}) \otimes A \to H^{\otimes r} \otimes A^{\otimes q}
   \end{split}
   \]
when $r \geq 1$.  When $r = 0$, the left and right $A$-action maps on $A^{\otimes q}$ are simply $\mu_{l,A}^q$ and $\mu_{r,A}^q$.  In particular, we can consider the Hochschild cochain complex $\Hoch^*(A, H^{\otimes r} \otimes A^{\otimes q})$.

\subsection{Comodule-bialgebra deformation tricomplex}
\label{subsec:def tri cb}

The \emph{deformation tricomplex} of the $H$-comodule-bialgebra $A$ is the first octant, cohomological tricomplex
   \[
   \Ccb^{***}(A) \,=\, \left\lbrace \Ccb^{p,\,q,\,r}(A), (-1)^{q+1}\dI^{p,\,q,\,r}, (-1)^{r+1}\dII^{p,\,q,\,r}, (-1)^{p+1} \dIII^{p,\,q,\,r}\right\rbrace
   \]
with
   \[
   \Ccb^{p,\,q,\,r}(A) \,=\,
   \begin{cases}
   0 & \text{ if $p = 0$ or $q = 0$}, \\
   \Bider(A) & \text{ if $(p,q,r) = (1,1,0)$}, \\
   \Hom(A^{\otimes p}, H^{\otimes r} \otimes A^{\otimes q}) & \text{ otherwise}.
   \end{cases}
   \]
The differential
   \[
   (-1)^{q+1}\dI^{p,\,q,\,r} \colon \Ccb^{p,\,q,\,r}(A) \to \Ccb^{p+1,\,q,\,r}(A)
   \]
is defined using the Hochschild coboundary
   \[
   \dI^{p,\,q,\,r} \,\buildrel \text{def} \over=\, \delta_h^p \colon
   \Hoch^p(A, H^{\otimes r} \otimes A^{\otimes q}) \to
   \Hoch^{p+1}(A, H^{\otimes r} \otimes A^{\otimes q}).
   \]

The differential
   \[
   (-1)^{r+1}\dII^{p,\,q,\,r} \colon \Ccb^{p,\,q,\,r}(A) \to \Ccb^{p,\,q+1,\,r}(A)
   \]
is defined using the alternating sum
   \[
   \dII^{p,\,q,\,r} \,=\, \sum_{i=0}^{q+1} \, (-1)^i \dII^{p,\,q,\,r} \lbrack i \rbrack,
   \]
where
   \[
   \dII^{p,\,q,\,r} \lbrack i \rbrack(\varphi) \,=\,
   \begin{cases}
   \left(\mu_{l,H}^r \otimes \mu_A^p \otimes \Id_{A^{\otimes q}}\right) \circ \left(\Id_H \otimes \tau_{(A^{\otimes p}, H^{\otimes r})} \otimes \Id_{A^{\otimes q}}\right) \circ \left(\rho^p \otimes \varphi\right) \circ \Delta_{A^{\otimes p}} & \text{ if $i = 0$}, \\
   \left(\Id_{H^{\otimes r} \otimes A^{\otimes (i-1)}} \otimes \Delta_A \otimes \Id_{A^{\otimes (q-i)}}\right) \circ \varphi & \text{ if $1 \leq i \leq q$}, \\
   \left(\mu_{r,H}^r \otimes \Id_{A^{\otimes q}} \otimes \mu_A^p\right) \circ \left(\Id_{H^{\otimes r}} \otimes \tau_{(A^{\otimes q}, H)} \otimes \Id_{A^{\otimes p}} \right) \circ \left(\varphi \otimes \rho^p\right) \circ \Delta_{A^{\otimes p}} & \text{ if $i = q+1$}.
   \end{cases}
   \]
Here $\rho^p \colon A^{\otimes p} \to H \otimes A^{\otimes p}$ and $\mu_A^p \colon A^{\otimes p} \to A$ are as defined in section \ref{subsec:ca} and \eqref{eq:mun}, respectively.  The map
   \[
   \Delta_{A^{\otimes p}} \colon A^{\otimes p} \to A^{\otimes p} \otimes A^{\otimes p}
   \]
is defined as
   \[
   \Delta_{A^{\otimes p}}(a_{1,p})
   \,=\, \sum \left(a_{1(1)} \otimes \cdots \otimes a_{p(1)}\right) \otimes \left(a_{1(2)} \otimes \cdots \otimes a_{p(2)}\right)
   \]
for $a_1, \ldots , a_p \in A$.

Finally, the differential
   \[
   (-1)^{p+1} \dIII^{p,\,q,\,r} \colon \Ccb^{p,\,q,\,r}(A) \to \Ccb^{p,\,q,\,r+1}(A)
   \]
is defined using the alternating sum
   \[
   \dIII^{p,\,q,\,r} \,=\, \sum_{i=0}^{r+1} \, (-1)^i \dIII^{p,\,q,\,r} \lbrack i \rbrack,
   \]
where
   \[
   \dIII^{p,\,q,\,r} \lbrack i \rbrack(\varphi) \,=\,
   \begin{cases}
   \left(\Id_{H^{\otimes r}} \otimes \rho^q\right) \circ \varphi & \text{ if $i = 0$}, \\
   \left(\Id_{H^{\otimes (r-i)}} \otimes \Delta_H \otimes \Id_{H^{\otimes (i-1)} \otimes A^{\otimes q}}\right) \circ \varphi & \text{ if $1 \leq i \leq r$}, \\
   \left(\Id_H \otimes \varphi\right) \circ \rho^p & \text{ if $i = r+1$}.
   \end{cases}
   \]

The six statements in Theorem \ref{thm:bicomplex mb} still hold in the comodule-bialgebra context, so $\Ccb^{***}(A)$ is indeed a tricomplex.

\subsection{Boundary planes}
\label{subsec:boundary plane cb}

As in the module-bialgebra case, the boundary planes of the deformation tricomplex $\Ccb^{***}(A)$  are either known or have been discussed in previous sections.  In fact:
\begin{enumerate}
\item The $p = 1$ plane $\Ccb^{1,\,*,\,*}(A)$ coincides with the deformation bicomplex $\Ccc^{**}(A)$ in which $A$ is regarded as an $H$-comodule-coalgebra (see section \ref{sec:cc}).
\item The $q = 1$ plane $\Ccb^{*,\,1,\,*}(A)$ coincides with the deformation bicomplex $\Cca^{**}(A)$ in which $A$ is regarded as an $H$-comodule-algebra (see section \ref{sec:ca}).
\item The $r = 0$ plane $\Ccb^{*,\,*,\,0}(A)$ coincides with the Gerstenhaber-Schack deformation bicomplex $\CGS^{**}(A)$ (denoted by $\widehat{C}_b^{\bullet,\, \bullet}(A,A)$ in \cite{gs1,gs2}) in which $A$ is regarded as only a bialgebra.
\end{enumerate}
Again, the only exception to the above remarks is the entry $\Ccb^{1,\,1,\,0}(A) = \Bider(A)$.

\subsection{Cup-products}
\label{subsec:cup cb}

Each cochain complex in the deformation tricomplex $\Ccb^{***}(A)$ is a DGA.

In the direction of $\dI$, fix $q \geq 1$ and $r \geq 0$.  Note that the module $H^{\otimes r} \otimes A^{\otimes q}$ is an associative algebra whose multiplication map is
   \[
   \mu_{H^{\otimes r} \otimes A^{\otimes q}} \,=\,
   \left(\mu_{H^{\otimes r}} \otimes \mu_{A^{\otimes q}}\right) \circ \left(\Id_{H^{\otimes r}} \otimes \tau_{(A^{\otimes q},\, H^{\otimes r})} \otimes \Id_{A^{\otimes q}}\right).
   \]
Regarding $H^{\otimes r} \otimes A^{\otimes q}$ as an $A$-bimodule, the three conditions in \eqref{eq:3conditions} hold.  It follows that $\Ccb^{*,\,q,\,r}(A) = (-1)^{q+1}\Hoch^*(A, H^{\otimes r} \otimes A^{\otimes q})$ admits a Hochschild $\cup$-product that makes it into a DGA.  This generalizes the $\cup$-product in $\Cca^{*,\,r}(A)$ discussed in section \ref{subsec:cup ca}.

In the direction of $\dII$, fix $p \geq 1$ and $r \geq 0$.  Then the cochain complex $\Ccb^{p,\,*,\,r}(A)$ is a DGA whose product is defined as
   \[
   \varphi \cup \psi \,\buildrel \text{def} \over= \,
   \left(\mu_{H^{\otimes r}} \otimes \Id_{A^{\otimes (q_1+q_2)}}\right) \circ
   \left(\Id_{H^{\otimes r}} \otimes \tau_{(A^{\otimes q_1},\, H^{\otimes r})} \otimes \Id_{A^{\otimes q_2}}\right) \circ
   (\varphi \otimes \psi) \circ \Delta_{A^{\otimes p}}
   \]
for $\varphi \in \Ccb^{p,\,q_1,\,r}(A)$ and $\psi \in \Ccb^{p,\,q_2,\,r}(A)$.  This generalizes the $\cup$-product in $\Ccc^{*,\,r}(A)$ discussed in section \ref{subsec:cup cc}.

In the direction of $\dIII$, fix $p \geq 1$ and $q \geq 1$.  Then the cochain complex $\Ccb^{p,\,q,\,*}(A)$ is a DGA whose product is defined as
   \[
   \varphi \cup \psi \,\buildrel \text{def} \over= \,
   \left(\Id_{H^{\otimes (r_1+r_2)}} \otimes \mu_{A^{\otimes q}}\right) \circ
   \left(\Id_{H^{\otimes r_1}} \otimes \tau_{(A^{\otimes q},\, H^{\otimes r_2})} \otimes \Id_{A^{\otimes q}}\right) \circ (\varphi \otimes \psi) \circ \Delta_{A^{\otimes p}}
   \]
for $\varphi \in \Ccb^{p,\,q,\,r_1}(A)$ and $\psi \in \Ccb^{p,\,q,\,r_2}(A)$.

\subsection{Total complex}
\label{subsec:tot cb}

The total complex $\Ccb^*(A)$, called the \emph{deformation complex of $A$}, of the deformation tricomplex $\Ccb^{***}(A)$ is defined as usual but with a shift of degree:
   \[
   \Ccb^n(A) \,\buildrel \text{def} \over=\, \bigoplus_{p+q+r \,=\, n+1} \Ccb^{p,\,q,\,r}(A).
   \]
In particular, we have that
   \[
   \begin{split}
   \Ccb^1(A) & \,=\, \Bider(A), \\
   \Ccb^2(A) & \,=\, \Hom(A, H \otimes A) \oplus \Hom(A^{\otimes 2}, A) \oplus \Hom(A, A^{\otimes 2}).
   \end{split}
   \]
The $n$th cohomology module of the deformation complex $\Ccb^*(A)$ is denoted by $\Hcb^n(A)$.

\subsection{Comodule-bialgebra deformation}
\label{subsec:def cb}

A \emph{deformation of $A$} as an $H$-comodule-bialgebra is a power series $\Theta = \sum_{n \geq 0} \theta_n t^n$ with $\theta_0 = (\rho, \mu_A, \Delta_A) \in \Ccb^2(A)$ and each $\theta_n = (\rho_n, \pi_n, \Delta_n) \in \Ccb^2(A)$, satisfying the five conditions for comodule-algebra deformation \eqref{eq:def ca} and comodule-coalgebra deformation \eqref{eq:def cc}.  The \emph{infinitesimal} of $\Theta$ is the linear coefficient $\theta_1$.

The notions of a \emph{formal automorphism of $A$} and \emph{equivalence} of deformations are defined exactly as in the module-bialgebra case, except that, in defining the latter, the conditions \eqref{eq:equiv ca} and \eqref{eq:equiv cc} are used.

The \emph{trivial deformation} is the deformation $\Theta = \theta_0 = (\rho, \mu_A, \Delta_A)$.  The $H$-comodule-bialgebra $A$ is said to be \emph{rigid} if every deformation of $A$ is equivalent to the trivial deformation.

The following result is the comodule-bialgebra analogue of Theorem \ref{thm:def ma}.

\pagebreak

\begin{thm}
\label{thm:def cb}
Let $A$ be an $H$-comodule-bialgebra.  Then the following statements hold.
\begin{enumerate}
\item The infinitesimal $\theta_1$ of any deformation $\Theta$ of $A$ is a $2$-cocycle in the deformation complex $\Ccb^2(A)$ whose cohomology class is determined by the equivalence class of $\Theta$.
\item If $\Hcb^2(A) = 0$, then $A$ is rigid.
\end{enumerate}
\end{thm}

\sgsp


\end{document}